# Stability and Boundedness of Solutions to Some Multidimensional Time-Varying Nonlinear Systems

Mark A. Pinsky

**Abstract**. Assessment of the degree of boundedness/stability of multidimensional nonlinear systems with time-dependent and nonperiodic coefficients is an important problem in various applied areas which has no adequate resolution yet. Most of the known techniques provide computationally intensive and conservative stability criteria in this field and frequently fail to estimate the regions of boundedness/stability of solutions to the corresponding systems. Recently, we outlined a new approach to this problem which is based on the analysis of solutions to a scalar auxiliary equation bounded from the above time histories of the norms of solutions to the original system. This paper develops a novel technique casting the auxiliary equation in a modified form which extends the application domain and reduces the computational hamper of our prior approach. Consequently, we developed more general boundedness/stability criteria and estimated trapping/stability regions for some multidimensional nonlinear systems with nonperiodic time – dependent coefficients that are common in various application domains. This let us to assess in target simulations the extent of boundedness/stability of some multidimensional, nonlinear and time - varying systems which were intractable with our prior technique.

**Key Words**: Time-varying nonlinear systems; nonperiodic coefficients; criteria of stability and boundedness of solutions; estimation of trapping/stability regions.

## 1. Introduction and Motivation Example

### 1.1 Background

Analysis of boundedness and stability of nonlinear systems with variable and nonperiodic coefficients plays a vital role in various engineering and natural science problems, which, for instance, are concerned with the design of controllers and observers. It appears that currently there are no confirmable necessary and sufficient conditions of local stability of the trivial solution to homogeneous systems of such kind, see, e.g., [1]-[7]. It was shown by A. M. Lyapunov [8] that under some additional conditions, the trivial solution to a homogeneous system with time-varying coefficients is asymptotically stable if the linearization of this system at zero is regular and its maximal Lyapunov exponent is negative, see also [1] for contemporary review of this subject. Yet, the confirmation of the former condition presents a considerable problem in applications. Subsequently, it was demonstrated by Perron [9] that arbitrary small perturbations can reverse sign of the Lyapunov exponents of a linearized system with time-varying coefficients and alter its stability. This attested that Lyapunov's regularity condition is essential. Consequent works were focused on examining the Lyapunov stability conditions through a review of the stability of Lyapunov exponents for linearized nonautonomous systems. Finally, the necessary and sufficient conditions of stability of Lyapunov exponents for solutions to the linearized system were given in [10] and [11]. However, it turns out that authentication of these conditions requires the apprehension of the fundamental set of solutions for the underlying systems, which is dubious in applications.

Application of the concept of generalized exponents, that was introduced in [6], provides sufficient conditions of stability of the trivial solution to time-varying nonlinear systems which, in principle, can be checked with the aid of numerical simulations. Still, the upper generalized exponent is larger than or equal to the maximal Lyapunov exponent, which heightens the conservatism of this more robust approach.

To retrieve a simple sufficient stability condition from the last approach, we firstly write the following equation with marked linearization,

$$\dot{x} = B(t)x + f_*(t,x) + F_*(t), \quad \forall t \geq t_0, \quad x(t,t_0,x_0) \in \mathbb{R}^n, \quad f_*(t,0) = 0$$
$$x(t_0,t_0,x_0) = x_0 \tag{1.1}$$

where functions, $f_* : [t_0,\infty) \times \mathbb{R}^n \to \mathbb{R}^n$, $F_* : [t_0,\infty) \to \mathbb{R}^n$ and matrix, $B : [t_0,\infty) \to \mathbb{R}^{n \times n}$ are continuous, $t_0 \in \mathrm{T} := [\zeta,\infty)$, $\zeta \in \mathbb{R}$, $x_0 \in H \subset \mathbb{R}^n$, $H$ is a neighborhood of $x \equiv 0$

---

Mark A. Pinsky, Department of Mathematics and Statistics, University of Nevada.Reno, Reno NV 89557, USA, e-mial:pinsky@unr.edu.

, $\mathbb{R}^n$ is a real $n$-dimensional Euclidean space, $F_*(t) = F_0 \eta(t)$, $\sup_{t \geq t_0} \|\eta(t)\| = 1$, $F_0 > 0$, $\|\cdot\|$ stands for induced 2-norm of a matrix or 2-norm of a vector, and $x(t,t_0,x_0):[t_0,\infty) \times T \times H \to \mathbb{R}^n$ is a solution to (1.1). To shorten the notation, we will write below that $x(t,t_0,x_0) \equiv x(t,x_0)$ and assume that (1.1) possesses a unique solution $\forall x_0 \in H$ and $\forall t \geq t_0$.

Due to continuity of the right side of (1.1), the last assumption implies that $x(t,x_0)$ is a continuous and continuously differential function that is bounded for any $t \in [t_0, t_*]$, $t_* < \infty$. Consequently, the reminder of this paper focuses on behavior of $x(t,x_0)$ for $t \to \infty$.

We also will assess stability of the trivial solution to the following homogeneous equation,
$$\dot{x} = B(t)x + f_*(t,x)$$
$$x(t_0, x_0) = x_0 \tag{1.2}$$
and its linear counterpart,
$$\dot{x} = B(t)x$$
$$x(t_0, x_0) = x_0 \tag{1.3}$$

Clearly, the solution to (1.3) can be written as $x(t,x_0) = W(t,t_0)x_0$, where $W(t,t_0) = w(t)w^{-1}(t_0)$ and $w(t)$ are transition and fundamental solution matrices for (1.3), respectively [1], [4].

Next, we write Lipschitz continuity condition for $f_*(t,x)$ as follows,
$$\|f_*(t,x)\| \leq l_1(t,t_0)\|x\|, \quad \forall x \in \Omega_1 \subset \mathbb{R}^n, \forall t \geq t_0, \tag{1.4}$$
where $\Omega_1$ is a bounded neighborhood of $x \equiv 0$ and $l_1(t,t_0) \leq \hat{l}_1 > 0$, $\forall t \geq t_0$ is a continuous function. In turn, let us assume that,
$$\|W(t,t_0)\| \leq Ne^{-v(t-t_0)}, \quad \forall t \geq t_0, \ N, v > 0, \tag{1.5}$$
and
$$N\hat{l}_1 - v < 0 \tag{1.6}$$
Then the trivial solution to (1.2) is asymptotically stable if (1.4), (1.5), and (1.6) hold [4, 6, 7].

A more general but less tractable condition on the norm of the transition matrices was presented in [1, 2] in the form,
$$\|W(t,t_0)\| \leq N \exp \int_{t_0}^{t} \Upsilon(s)ds, \ \forall t > t_0, \ \Upsilon:[t_0,\infty) \to \mathbb{R}, \ N > 0 \tag{1.7}$$
where it was shown that (1.4), (1.7), and the following condition,
$$\limsup_{t \to \infty}(1/(t-t_0))\int_{t_0}^{t} \Upsilon(s)ds + N\hat{l}_1 < 0 \tag{1.8}$$
embracing the asymptotic stability of the trivial solution to (1.2). More compelling stability conditions of (1.2) were given in [12], yet verification of these conditions in applications can be challenging as well.

In the control literature, the analysis of stability of nonautonomous nonlinear systems was frequently aided by application of the Lyapunov function method [15 – 23]. Nonetheless, the application of this methodology is strenuous in this area since adequate Lyapunov functions are rarely known for multidimensional and nonlinear systems with nonperiodic coefficients.

The problem of estimating the stability regions of autonomous nonlinear systems has attracted numerous publications in the last few decades [24]-[38], but these techniques fail for systems with time-varying coefficients. To our knowledge, the problem of estimating the trapping regions for nonautonomous nonlinear systems has not been virtually addressed in the current literature.

Assessment of stability of nonlinear systems based on the analysis of convergence of adjacent trajectories was developed in [39]. However, to our knowledge, this approach rarely leads to computationally sound stability conditions for multidimensional systems.

In [13], we developed a novel technique for estimation of upper bounds of the norms of solutions to (1.1) or (1.2) and use it to distinguish the boundedness/stability criteria and estimate the trapping/stability regions for these systems. Consequently, under normalizing condition, $\|w(t_0)\|=1$, we derive the following inequality,

$\|x(t,x_0)\| \leq X(t,X_0)$, $\forall t \geq t_0$, where $x(t,x_0)$ is a solution to (1.1) and $X(t,X_0):[t_0,\infty)\times\mathbb{R}_{\geq 0} \to \mathbb{R}_{\geq 0}$ ($\mathbb{R}_{\geq 0}$ is a set of nonnegative real numbers) is a solution to a scalar equation we derived in [13],

$$\dot{X} = p(t)X + c(t)\|f_*(t,x(t,x_0)) + F_*(t)\|$$
$$X(t_0,X_0) = \|w^{-1}(t_0)x_0\| = X_0 \tag{1.9}$$

Note that,
$$p(t) = d\left(\ln\|w(t)\|\right)/dt \tag{1.10}$$

and
$$c(t) = \|w(t)\|\|w^{-1}(t)\| = \sigma_{\max}(w)/\sigma_{\min}(w) \tag{1.11}$$

In turn, [13] has also introduced a nonlinear extension of the Lipschitz continuity condition which was presented as follows,

$$\|f_*(t,x)\| \leq L(t,\|x\|), \quad \forall x \in \Omega_2^x \subset \mathbb{R}^n, \quad \forall t \geq t_0 \tag{1.12}$$

where $L:=[t_0,\infty)\times\mathbb{R}_{\geq 0} \to \mathbb{R}_{\geq 0}$ is continuous function in $t$ and $\|x\|$, and locally Lipschitz in $\|x\|$, $L(t,0)=0$, $\Omega_2^x$ is a bounded neighborhood of $x \equiv 0$. Note that [13] defines (1.12) in a closed – form if $f_*$ is either a piece-wise polynomial function in $x$ or can be approximated by such function with, e.g., bounded in $\Omega_2^x$ error term, which, for instance, can be written in Lagrange form. In the former case, $\Omega_2^x \equiv \mathbb{R}^n$ and in the latter case, $\Omega_2^x \equiv \mathbb{R}^n$ if the error term is also globally bounded for $\forall x \in \mathbb{R}^n$. Thus, condition, $\Omega_2^x \equiv \mathbb{R}^n$ holds for a large set of nonlinear systems emerging in science and engineering applications and we are going to use it in the reminder of this paper.

Let us illustrate how to define $L(t,\|x\|)$ for a simple but representative vector function. Assume that $x = [x_1 \ x_2]^T$ and a vector-function $f \in \mathbb{R}^2$ is defined, e.g., as follows, $f = \begin{bmatrix} a_1(t)x_1^3 x_2^2 & a_2(t)x_2^2 \end{bmatrix}^T$. Then, $\|f\|_2 \leq \|f\|_1 \leq |a_1||x_1^3||x_2^2| + |a_2||x_2^2| \leq |a_1|\|x\|^5 + |a_2|\|x\|^2$, where we use that $|x_i^n| \leq \|x\|^n$, $i=1,2$, $n \in \mathbb{N}$ and $\mathbb{N}$ is a set of positive integers. Note that additional and more complex examples of this kind are provided in Section 5, see also [13] and [14].

Consequently, using (1.12), we reviewed (1.9) in the following more tractable form,
$$\dot{z} = p(t)z + c(t)\left(L(t,z) + \|F_*(t)\|\right)$$
$$z(t_0,z_0) = z_0 = \|w^{-1}(t_0)x_0\| \tag{1.13}$$

To reduce the notation, we adopt in (1.13) and throughout this paper that $z(t,t_0,z_0) = z(t,z_0)$.

In turn, it was shown in [13] that
$$\|x(t,x_0)\| \leq X(t,X_0) \leq z(t,z_0) \tag{1.14}$$

where $x(t,x_0)$, $X(t,X_0)$ and $z(t,z_0)$ are solutions of (1.1), (1.9) and (1.13), respectively.

Thus, (1.14) can be used to bound from above the norm of solutions of (1.1) by matching solutions of the scalar equation (1.13).

It turned out that (1.13) provides sound estimates of the trapping/stability regions under the assumption that only the bound of $\|f_*(t,x)\|$ is known – a frequent pronouncement in the control literature [3]-[5]. Nonetheless, such estimates may become more conservative if $f_*(t,x)$ is defined explicitly. Consequently, we refined this methodology in [14], where (1.13) was used to estimate the error of successive approximations of solutions to (1.1) or (1.2) stemming from the trapping/stability regions of these systems. This modified approach enhanced our boundedness/stability criteria and delivered approximations that increased successively the accuracy of estimation of the boundaries of the corresponding trapping/stability regions.

Nonetheless, the methodologies developed in [13] and [14] works under the condition that $c(t) < \infty$, $\forall t \geq t_0$, which considerably limits the scope of its applications. Nonetheless, it follows from (1.11) that frequently, $\lim_{t \to \infty} c(t) = \infty$ even if $A$ is a stable and time-invariant matrix.

The current paper lifts this limitation for a practically important class of nonlinear systems with time-varying and nonperiodic coefficients. Its main contribution is in the development of a modified auxiliary equation with $c(t) \equiv 1$ and $p(t) \equiv const$ under some conditions that are frequently met in various applications. This new technique prompts the development of novel criteria of bondedness/stability and estimation of bondedness/stability regions for a wide class of systems that were intractable to our former methodology [13] and voids elaborate simulations of $w(t)$, $p(t)$ and $c(t)$ that were required previously.

To make this paper more inclusive, we present the definitions of some standard principles of stability theory which are going to be used in the reminder of this paper. Note that the standard definitions of Lyapunov stability and asymptotic stability for time-varying nonlinear systems that are accepted below can be found, e.g., in [1]-[5]. In the reminder of this paper, we will call these properties shortly either stability or asymptotic stability.

We begin with a formal definition of Lyapunov exponents, where we adopt the exposition of these quantities made in [1] and [5].

For equation (1.2), the Lyapunov exponents $\phi_i(x(t,x_0)) = \lim_{t_0 < t \to \infty} \sup (t-t_0)^{-1} \ln \|x(t,x_0)\|$, $i \leq n$, where $x(t,x_0)$ is a solution to (1.2). The maximal Lyapunov exponent, $\phi = \max_i \{\phi_i(x(t,x_0))\}$ measures the maximal rate of exponential growth/decay of the corresponding solutions as $t \to \infty$ and bears a pivoting role in stability theory. For a linear system (1.3), the Lyapunov exponents are defined as, $\phi_{l,i} = \lim_{t_0 < t \to \infty} \sup (t-t_0)^{-1} \ln \sigma_i(\|w(t)\|)$, $i \leq n$, where $\sigma_i$ are the singular values of the fundamental solution matrix of (1.3); and for linear systems, $\max_i \phi_{l,i} = \phi$ [1]. This let to bound the transition matrix of (1.3) as follows,

$$\|W(t,t_0)\| \leq D(\vartheta) e^{-\chi(t-t_0)}, \ \forall t \geq t_0, \ D > 0, \ -\chi = \phi + \vartheta \qquad (1.15)$$

where $\vartheta$ is a small positive number [6].

In order, let us bring the definition of the comparison principle [4] which is frequently used below. Consider a scalar differential equation

$$\dot{u}_1 = g(t,u_1), \ \forall t \geq t_0, \ u_1(t,u_{10}) \in \mathbb{R}$$
$$u_1(t_0,u_{10}) = u_{10} \qquad (1.16)$$

where function $g(t,u_1)$ is continuous in $t$ and locally Lipschitz in $u_1$ for $\forall u_1 \in \wp \subset \mathbb{R}$. Suppose that the solution to this equation, $u_1(t,u_{10}) \in \wp$, $\forall t \geq t_0$. Next, consider a differential inequality,

$$D^+u_2 \leq g(t, u_2), \ \forall t \geq t_0$$
$$u_2(t_0, u_{20}) = u_{20} \leq u_{10}$$

where $D^+u_2$ denotes the upper right-hand derivative in $t$ of $u_2(t, u_{20})$ [4] and $u_2(t, u_{20}) \in \wp, \ \forall t \geq t_0$. Then, $u_1(t, u_{10}) \geq u_2(t, u_{20}), \ \forall t \geq t_0$.

Now we present for convenience some conventional definitions of the trapping/stability regions as follows.

**Definition 1**. A connected and compact set of all initial vectors, $\mathfrak{I}_1(t_0)$ that includes zero-vector, is called a trapping region of equation (1.1) about the origin of coordinates frame $x \equiv 0$ if condition $x_0 \in \mathfrak{I}_1(t_0)$ implies that $x(t, x_0) \in \mathfrak{I}_1(t_0), \ \forall t > t_0$.

Clearly, this definition acknowledges that $\mathfrak{I}_1$ is the invariant set of (1.1) containing zero.

**Definition 2**. A connected and open set of all initial vectors, $\mathfrak{I}_2(t_0)$, that includes zero-vector, is called a region of stability of the trivial solution to (1.2) if condition $x_0 \in \mathfrak{I}_2(t_0)$ implies that $x(t, x_0)$ is stable.

**Definition 3**. A connected and open set of all initial vectors, $\mathfrak{I}_3(t_0)$, that includes zero-vector, is called a region of asymptotic stability of the trivial solution to (1.2) if condition, $x_0 \in \mathfrak{I}_3(t_0)$ implies that $\lim_{t \to \infty} x(t, x_0) = 0$.

### 1.2 Motivation Example

The first part of this section derives a novel auxiliary equation with $c(t) = 1$ and $p(t) = const$ for a simple planar system. The second part – inference stability and estimate the stability region of the derived auxiliary equation and extends these inferences to the stability assessment of our planar model (1.2).

Let us firstly apply our former methodology [13] to a simple case, where $B = diag(\lambda_1 \ \lambda_2), \lambda_{1,2} \in \mathbb{R}, \lambda_1 > \lambda_2$. Obviously, in this case $w(t) = diag\left(e^{\lambda_1(t-t_0)} \ e^{\lambda_2(t-t_0)}\right)$, $w^{-1}(t) = diag\left(e^{-\lambda_1(t-t_0)} \ e^{-\lambda_2(t-t_0)}\right)$ which implies that $\|w(t)\| = e^{\lambda_1(t-t_0)}$, $\|w^{-1}(t)\| = e^{\lambda_2(t-t_0)}$. Thus, $c(t) = \|w(t)\|\|w^{-1}(t)\| = e^{(\lambda_1-\lambda_2)(t-t_0)}$ and $\lim_{t \to \infty} c(t) = \infty$ if $\lambda_1 \neq \lambda_2$ which makes our prior estimates of $\|x(t, x_0)\|$, which is based on (1.13), overconservative for large values of $t$. Nonetheless, $c(t) \equiv 1$ if $\lambda_{1,2}$ are complex conjugate and our prior inferences hold in this case.

Let us break out our current approach into a sequence of straightforward steps for a simple version of system (1.1) where $x \in \mathbb{R}^2$, $B(t) = A + G_*(t)$, $A \in \mathbb{R}^{2 \times 2}$, $A = const$, $A \neq 0$, $G_*(t) \in \mathbb{R}^{2 \times 2}$ and $\lim_{t \to \infty} (t - t_0)^{-1} \int_{t_0}^{t} G_*(s) ds = 0$. Assume also that matrix $A$ is diagonalizable and eigenvalue and eigenvector matrices of $A$ are $\Lambda = diag(\lambda, \lambda^*)$, $V = [v_{kj}]$, $k, j = 1, 2$, where $\lambda = \alpha + i\beta$, $\lambda^* = \alpha - i\beta$, $f_* = \begin{pmatrix} 0 & \varepsilon x_2^3 \end{pmatrix}^T$, $F_* = \begin{pmatrix} 0 & F_{*2} \end{pmatrix}^T$.

Equation (1.13) can be applied to the stability analysis of such a planar system if (1.11) implies that $c(t) < \infty, \forall t \geq t_0$ which, in general, is difficult to warranty in advance. Thus, we outline a different technique recasting the auxiliary equation in the form, where $c(t) \equiv 1$ and $p(t) \equiv const$.

1. Firstly, let us write our planar model of system (1.1) in eigenbasis of $A$ as follows

$$\dot{y} = \Lambda y + G(t) y + V^{-1} \left( 0 \quad -\varepsilon \left( \sum_{j=1}^{2} v_{2j} y_j \right)^3 \right)^T + F_2(t)$$

where $y \in \mathbb{C}^2$, $\mathbb{C}^2$ is a two-dimensional space of complex numbers, $x = Vy$, $G(t) = V^{-1} G_*(t) V$ and $F(t) = V^{-1} F_*(t)$.

2. Next, we rewrite the last equation as follows,

$$\dot{y} = (\lambda I + ib) y + (a - \lambda I) y + G(t) y + V^{-1} \left( 0 \quad -\varepsilon \left( \sum_{j=1}^{2} v_{2j} y_j \right)^3 \right)^T + F(t) \quad (1.17)$$

where $I$ is an identity matrix, $b = diag(\beta \quad -\beta)$, $a = diag(\alpha \quad \alpha)$ and $\lambda \in \mathbb{R}$ is going to be determine subsequently.

In order, let us select a linear subsystem of (1.17) with an underlined diagonal matrix as,

$$\dot{y} = (\lambda I + ib) y \quad (1.18)$$

For this subsystem, a diagonal fundamental solution matrix can be written as follows,
$w(t) = e^{(\lambda I + ib)(t - t_0)} = e^{\lambda(t - t_0)} e^{ib(t - t_0)}$ which implies that $\|w(t)\| = e^{\lambda(t - t_0)} \|e^{ib(t - t_0)}\| = e^{\lambda(t - t_0)}$ since $\|e^{ib(t - t_0)}\| = 1$.
Same reasoning yields that, $\|w^{-1}(t)\| = e^{-\lambda(t - t_0)}$ which shows that $c(t) = \|w(t)\| \|w^{-1}(t)\| = 1$, and, due to (1.10), $p(t) = \lambda = const$.

Thus, if we use (1.18) as the underlying linear system, then the first addition in right side of the auxiliary equation (1.13), which is derived for (1.17), is $p(t) z = \lambda z$.

3. Further application of (1.13) to (1.17) brings the second term in the right side of our modified auxiliary equation in the following form,

$$\left\| (a - \lambda I) y + G(t) y + V^{-1} \left( 0 \quad -\varepsilon \left( \sum_{j=1}^{2} v_{2j} y_j \right)^3 \right)^T + F(t) \right\| \leq$$

$$\left( \|a - \lambda I\| + \|G(t)\| \right) \|y\| + \|V^{-1}\| \left\| \left( 0 \quad -\varepsilon \left( \sum_{j=1}^{2} v_{2j} y_j \right)^3 \right)^T \right\| + \|F(t)\|$$

Then, the utility of standard inequality yields

$$\left\| \left( 0 \quad -\varepsilon \left( \sum_{j=1}^{2} v_{2j} y_j \right)^3 \right)^T \right\| \leq abs\left( \varepsilon \left( \sum_{j=1}^{2} v_{2j} y_j \right)^3 \right) \leq \varepsilon \left( abs \sum_{j=1}^{2} v_{2j} y_j \right)^3 \leq \varepsilon \left( \sum_{j=1}^{2} abs(v_{2j}) abs(y_j) \right)^3 \leq$$

$$\varepsilon \left( \sum_{j=1}^{2} abs(v_{2j}) \|y\| \right)^3 \leq \varepsilon \left( \sum_{j=1}^{2} abs(v_{2j}) \right)^3 \|y\|^3$$

where $abs(q) = \sqrt{a^2 + b^2}$ if $q = a + ib$, and $abs(q) = |q|$ if $\text{Im}(q) = 0$.

This let us to write (1.13) as follows,

$$\dot{z} = (\lambda + \|a - \lambda I\|) z + G(t) z + \sigma z^3 + \|F(t)\|$$
$$z(0) = \|V^{-1} x_0\|$$
(1.19)

where $\sigma = \varepsilon \|V^{-1}\| \left( \sum_{j=1}^{2} abs(v_{2j}) \right)^3$.

4. Lastly, we select $\lambda$ to maximize the degree of stability of a linear equation $\dot{z} = (\lambda + \|a - \lambda I\|) z$ which minimizes conservatism of our estimates and reduces to the following condition,

$$\min_{\lambda} (\lambda + \|a - \lambda I\|) = \min_{\lambda} (\lambda + |\alpha - \lambda|)$$
(1.20)

since $a$ is a diagonal matrix with equal eigenvalues. Resolving (1.20) yields that $\lambda = \alpha$ which casts (1.19) in the following form,

$$\dot{z} = (\alpha + \|G(t)\|) z + \sigma z^3 + \|F(t)\|, \ \forall t \geq t_0$$
$$z(t_0, z_0) = z_0 = \|V^{-1} x_0\|$$
(1.21)

For further references, we write a homogeneous counterpart to (1.21) as follows,

$$\dot{z} = (\alpha + \|G(t)\|) z + \sigma z^3, \ \forall t \geq t_0$$
$$z_0 = \|V^{-1} x_0\|$$
(1.22)

Let us recall now that, due to (1.14), $\|x(t, x_0)\| \leq z(t, z_0)$, $z_0 = \|V^{-1} x_0\|$, $\forall t \geq t_0$, where $x(t, x_0)$ and $z(t, z_0)$ are solutions to planar versions of either equation (1.1) or (1.2) and equation (1.21) or (1.22), respectively.

To further abridge the stability analysis of the trivial solution to (1.22), we are going to develop for this equation its linear, scalar and, thus, integrable upper bound. Due to the comparison principle [4], solutions to such equations resolve the essential boundedness/stability properties of equation (1.22). In general, Lipschitz continuity condition (1.4) can be used to develop such a linear equation. Nonetheless, in our case, we apply a less conservative bound, $z^3 \leq (\hat{z})^2 z$, $z \in [0, \hat{z}]$, $\hat{z} > 0$. Application of the last inequality brings the following linear equation

$$\dot{Z} = s(t, \hat{z}) Z, \ \forall t \geq t_0$$
$$Z(t_0, z_0, \hat{z}) = z_0 = \|V^{-1} x_0\|$$
(1.23)

where $s(t, \hat{z}) = \alpha + \|G(t)\| + \sigma (\hat{z})^2$. Note that solutions to (1.23) can be used to estimate solutions to (1.22) if $z \in [0, \hat{z}]$.

In turn, to simplify the assessment of stability of (1.23), we set that $\sup_{t_0 \leq t} s(t, \hat{z}) = -\alpha_* + (\hat{z})^2 \sigma + G_s = -\overline{s}(\hat{z})$, where we assumed that $\alpha = -\alpha_*$, $\alpha_* \geq 0$ and $G_s = \sup_{\forall t \geq t_0} \|G(t)\| < \infty$.

Subsequently, the solutions to (1.23) can be bounded as follows,

$$Z = z_0 \exp\left( \int_{t_0}^{t} s(\tau, \hat{z}) d\tau \right) \leq z_0 \exp\left( -\int_{t_0}^{t} \overline{s}(\hat{z}) d\tau \right) \leq z_0 \exp\left( -\overline{s}(\hat{z})(t - t_0) \right)$$ which implies that (1.22) is asymptotically stable if

$$\overline{s}(\hat{z}) > 0, \quad \hat{z} \in (0, \hat{z}_{max}) \tag{1.24}$$

where $\hat{z}_{max}$ is the maximal value of $\hat{z}$ for which (1.24) holds.

The last assessment can be drawn through the utility of a more general reasoning which is repeatedly applied in this paper. In fact, determine a comparison equation to (1.23) as follows.

$$\dot{\overline{Z}} = -\overline{s}(\hat{z})\overline{Z}, \; \forall t \geq t_0$$
$$\overline{Z}(t_0, z_0, \hat{z}) = z_0 = \|V^{-1}x_0\| \tag{1.25}$$

Clearly, the right side of (1.25) bounds from above the right side of (1.23), whereas the initial values for both equations are equal. Thus, $Z(t, z_0, \hat{z}) \leq \overline{Z}(t, z_0, \hat{z}) = z_0 \exp(-\overline{s}(\hat{z})(t - t_0))$, $\forall t \geq t_0$ due to comparison principle.

Next, we inference that (1.24) implies the asymptotic stability of the trivial solution to (1.22). In fact, due to the comparison principle,

$$z(t, z_0) \leq Z(t, z_0, \hat{z}), \; \forall t \geq t_0 \tag{1.26}$$

if $z(t, z_0) \in [0, \hat{z}]$, $\forall t \geq t_0$, where $z(t, z_0)$ and $Z(t, z_0, \hat{z})$ are solutions to (1.22) and (1.23).

Let us show that $z(t, z_0) \in [0, \hat{z})$, $\forall t \geq t_0$ if $z_0 < \hat{z}$ and (1.24) holds. Clearly, $z(t_0, z_0) = z_0 \in [0, \hat{z})$ under our assumption. Assume, on the contrary, that $t = t_* > t_0$ is the smallest value of $t$ for which $z(t_*, z_0) = \hat{z}$. Then, (1.26) holds for $\forall t \in [t_0, t_*]$ and, thus,

$z(t_*, z_0) \leq Z(t_*, z_0, \hat{z}) \leq \overline{Z}(t_*, z_0, \hat{z}) = z_0 \exp(-\overline{s}(\hat{z})(t_* - t_0)) < z_0 < \hat{z}$. Hence, $z(t, z_0) < \hat{z}$, $\forall t \geq t_0$ if (1.24) holds and $z_0 < \hat{z}$. Consequently, the trivial solutions to (1.22) and, in turn, to our planar model of equation (1.2) are asymptotically stable.

Furthermore, solving for $\hat{z}$ equation $\overline{s}(\hat{z}) = 0$ trumps estimation of the range of values of $\hat{z}$ which assure asymptotic stability of (1.25); i.e., values of $\hat{z} \in [0, \hat{z}_{max})$, where

$$\hat{z}_{max} = ((\alpha_* - G_s)/\sigma)^{1/2}$$

The last formula leads to the estimation of the stability regions for (1.22) and the corresponding model of equation (1.2). Really, (1.24) and, subsequently, (1.26) hold if $z_0 < \hat{z} \in (0, \hat{z}_{max})$. Hence, the last inequality implies that (1.25) and, in turn, (1.23) are asymptotically stable under that condition which, in turn, yields that the trivial solution to (1.22) is asymptotically stable as well, and a solution to (1.22) approaches zero as $t \to \infty$ if $z_0 \in (0, \hat{z}_{max})$. In order, this condition ensures asymptotic stability of the trivial solutions to our simplified model of (1.2) which relates to (1.22). Hence, due to the second relation in (1.22), the region of asymptotic stability of the trivial solution to our planar model of (1.2) includes interior points of the ellipsoid, which is defined as follows, $\|V^{-1}x_0\| = \hat{z}_{max}$.

Section 3 generalizes this last formula and derives the boundedness/stability criteria and estimates the trapping/stability regions for multidimensional equations (1.1) and (1.2).

## 2. Modified Auxiliary Equation

This section derives of a modified auxiliary equation with $c(t) = 1$ and $p(t) = const$ for a broad class of nonlinear systems that are frequently found in applications. Solutions to the auxiliary equation bound from above the matching solutions to (1.1) or (1.2) under some conditions which we specify below in a theorem which encapsulates our inferences and the underlined assumptions.

1. At this point, we assume that the average of $B(t)$ exists and is defined as follows

$$A(t_0) = \limsup_{t \to \infty, \forall t \geq t_0} (t-t_0)^{-1} \int_{t_0}^{t} B(s) ds < \infty, \ \forall t_0 \in T.$$ Additionally, we assume that $A(t_0) \in \mathbb{R}^{n \times n}$ is nonzero

and diagonalizable matrix $\forall t_0 \in T$ except possibly for some isolated values of $t_0$.

This last condition is met for any generic set of matrices $A(t_0)$ which depends upon a scalar parameter, i.e., such a set of matrices that is structurally stable under small perturbations, see, e.g., [41], chapter 6, pp.235-256 for more details.

In order, we set that complex conjugate eigenvalues of $A(t_0)$, $\lambda_k(t_0) = \alpha_k(t_0) \pm i\beta_k(t_0)$, $i = \sqrt{-1}$, $1 \leq k \leq n_1 \leq n$, $n_1 \in \mathbb{N}$, and real eigenvalues of $A(t_0)$, $\lambda_k(t_0) = \alpha_k(t_0)$, $n_1 < k \leq n$, where $\alpha_k \in \mathbb{R}$ and $\beta_k \in \mathbb{R}_{>0}$, where $\mathbb{R}_{>0}$ is a set of positive numbers. Additionally, we presume that $\alpha_k \geq \alpha_{k+1}$, $k \in [1, n-1]$ and also define a square diagonal matrix, $\Lambda = \alpha + i\beta$ with $\alpha = diag(\alpha_1, \alpha_1, ..., \alpha_{n_1}, \alpha_{n_1}, \alpha_{n_1+1}, ..., \alpha_n)$, $\beta = diag(\beta_1, -\beta_1, ..., \beta_{n_1}, -\beta_{n_1}, 0, ..., 0)$, $\alpha, \beta \in \mathbb{R}^{n \times n}$

2. Next, we write (1.1) as follows,

$$\dot{x} = Ax + G_*(t)x + f_*(t,x) + F_*(t), \ \forall t \geq t_0$$
$$x(t_0) = x_0 \tag{2.1}$$

where $G_*(t) = B(t) - A(t_0)$ is a zero-mean and continuous matrix.

Note that the immediate application of (1.13) to (2.1), which is based on utility of the fundamental matrix of solutions to equation, $\dot{x} = Ax$, i.e., $w(t) = \exp(At)$, frequently fails since in this case $c(t) = \exp((\alpha_1 - \alpha_n)t)$ and $\lim_{t \to \infty} c(t) = \infty$ if $\alpha_1 \neq \alpha_n$.

To escape this shortcoming, we rewrite (2.1) into the eigenbasis of $A$ as follows,

$$\dot{y} = \Lambda y + G(t)y + f(t,y) + F(t), \ \forall t \geq t_0, \ y \in \mathbb{C}^n$$
$$y(t_0, y_0) = y_0 = V^{-1}x_0 \tag{2.2}$$

where $y = V^{-1}x$, $V \in \mathbb{C}^{n \times n}$ is eigenmatrix of $A$, $\mathbb{C}^{n \times n}$ is $n^2$- dimensional space of complex numbers, $G = V^{-1}G_*V$, $f(t,y) = V^{-1}f_*(t,Vy)$ and $F(t) = V^{-1}F_*$. To reduce the notation, we adopted in (2.2) that $y(t, t_0, y_0) = y(t, y_0)$.

Subsequently, we rewrite (2.2) as,

$$\dot{y} = (\lambda I + i\beta)y + (\alpha - \lambda I)y + G(t)y + f(t,y) + F(t)$$
$$y_0 = V^{-1}x_0 \tag{2.3}$$

where $\lambda \in \mathbb{R}$ is defined below. Then, we select, $\dot{y} = (\lambda I + i\beta)y$ as our underlying linear equation with diagonal matrix, $\lambda I + i\beta = const$ and write the fundamental matrix of solutions to this equation as follows,

$$w(t) = \exp(\lambda I + i\beta)(t - t_0) = e^{\lambda(t-t_0)} \exp(i\beta)(t - t_0) \tag{2.4}$$

where $I$ is identity matrix. Henceforth, $\|w(t)\| = e^{\lambda(t-t_0)}$, $\|w^{-1}(t)\| = e^{-\lambda(t-t_0)}$ since $\|\exp i\beta(t-t_0)\| = 1$. In turn, application of (1.10) and (1.11) yields that $p(w(t)) = \lambda$ and $c(w(t)) = 1$ if $w(t)$ is defined by (2.4).

3. Next, adopting (2.4), we write equation (1.13) for (2.3) as follows,

$$\dot{z} = \left(\lambda + \|\alpha - \lambda I\|\right) z + \|G(t)\| z + L(t, z) + \|F(t)\|, \quad \forall \lambda \in \mathbb{R}$$
$$z(t_0, z_0) = z_0 = \|V^{-1} x_0\|$$
(2.5)

where $z(t, z_0) \geq \|y(t, y_0)\|$, $\forall t \geq t_0$, $y \in \Omega_2^y \subset \mathbb{C}^n$, $\Omega_2^y$ is a bounded neighborhood of $y \equiv 0$, $L := [t_0, \infty) \times \mathbb{R}_{\geq 0} \to \mathbb{R}_{\geq 0}$ is continuous in $t$ and $z$ and is locally Lipschitz in $z$ function with $L(t, 0) \equiv 0$. Note that this function can be developed through application of inequality (1.12) to function, $f(t, y) = V^{-1} f_*(t, Vy)$. Hence, (1.12) in this case takes the following form,

$$\|f(t, y)\| \leq L(t, \|y\|), \quad \forall y \in \Omega_2^y, \quad \forall t \geq t_0$$
(2.6)

On developing of $L(t, \|y\|)$ for polynomial vector – fields see examples in Sections 1.1 and 1.2 as well and more inclusive examples in Section 5 as well as examples in [13] and [14]. As we mentioned prior, inequality (2.6), for instance, holds if $f(t, y)$ is a piece – wise polynomial function or can be approximated by such function with bounded in $\Omega_2^y$ Lagrange-type error term. In the former case, $\Omega_2^y \equiv \mathbb{C}^n$ and in the latter case $\Omega_2^y \equiv \mathbb{C}^n$ if the error term is bounded in $\mathbb{C}^n$.

For straightforward further references, we assume in reminder of this paper that $\Omega_2^y \equiv \mathbb{C}^n$.

4. In the sequel, we define $\lambda$ using the following condition, $\min_\lambda \left(\lambda + \|\alpha - \lambda I\|\right)$ which maximizes the degree of stability of a scalar linear equation, $\dot{z} = \left(\lambda + \|\alpha - \lambda I\|\right) z$. Since $\alpha - \lambda I$ is a diagonal matrix, $\|\alpha - \lambda I\| = \max_k |\alpha_k - \lambda|$ which yields that,

$$\min_\lambda \left(\lambda + \|\alpha - \lambda I\|\right) = \min_\lambda \left(\lambda + \max_k |\alpha_k - \lambda|\right)$$

Clearly, the last condition yields that $\lambda = \alpha_n$ and, consequently, $\left(\lambda + \|\alpha - \lambda I\|\right)_{|\lambda = \alpha_n} = \alpha_1$, see proof in appendix of this paper.

This let us to write (2.5) as follows,

$$\dot{z} = \left(\alpha_1 + \|G(t)\|\right) z + L(t, z) + \|F(t)\|$$
$$z(t_0, z_0) = z_0 = \|V^{-1} x_0\|$$
(2.7)

5. To develop a less conservative version of (2.7), we set that matrix, $D(t) = \mathrm{Im}(diag G(t))$, $D \in \mathbb{R}^{n \times n}$, $\beta_+(t) = \beta + D(t)$ and $G_- = G - iD$, and rewrite (2.3) as follows,

$$\dot{y} = \left(\lambda I + i\beta_+(t)\right) y + \left(\alpha - \lambda I\right) y + G_-(t) y + f(t, y) + F(t)$$
$$y(t_0, y_0) = V^{-1} x_0$$

The fundamental matrix of solution for equation, $\dot{y} = \left(\lambda I + i\beta_+(t)\right) y$ again can be written as follows, $w_+(t) = \exp\left(\lambda I + i\beta_+(t)\right)(t - t_0)$ since $\beta_+$ is a diagonal matrix. As prior, this implies that $\|w_+(t)\| = e^{\lambda(t - t_0)}$, $\|w_+^{-1}(t)\| = e^{-\lambda(t - t_0)}$, $c(w_+(t)) = 1$ and $p(w_+(t)) = \alpha_1$. Lastly, a less conservative counterpart of (2.7) can be written as follows,

$$\dot{z} = \left(\alpha_1 + \|G_-(t)\|\right) z + L(t, z) + \|F(t)\|$$
$$z(t_0, z_0) = z_0 = \|V^{-1} x_0\|$$
(2.8)

In order, we encapsulate our derivations and the underlying assumptions in the following

**Theorem 1**. Assume that $B(t)$ is a continuous matrix, the average of $B(t)$, i.e., matrix $A(t_0) \neq 0$ exists and is diagonalizable for $\forall t_0 \in T$ except possibly some isolated values of $t_0$. Also assume that function $F_*(t)$ is continuous and function $f_*(t,x)$ is continuous in both variables and inequality (2.6) holds with $\Omega_2^y \equiv \mathbb{C}^n$. Lastly, we assume that equation (1.1), (2.7) and (2.8) possess unique solutions for $\forall t \geq t_0$, $\forall x_0 \in H \subset \mathbb{R}^n$, $\forall z_0 = \|V^{-1}x_0\|$. Then, the norms of solutions to equation (1.1) are bounded by matching solutions to the scalar equations (2.7) or (2.8) as follows,

$$\|x(t,x_0)\| \leq \|V\| z(t,z_0), \; z_0 = \|V^{-1}x_0\|, \; \forall t \geq t_0, \; x_0 \in H \qquad (2.9)$$

**Proof**. The proof of this statement directly follows from the above assumptions and inferences made prior in this section. In fact, our first set of conditions allows to set up equations (2.1) – (2.4). The next set of conditions implies inequality (2.6) with $\Omega_2^y \equiv \mathbb{R}^n$ and, in turn, equations (2.5), (2.7), and (2.8). The last assumption assures that functions $x(t,x_0)$ and $z(t,z_0)$ exist for $\forall t \geq t_0$ and the specified initial values □

Note that the above inferences naturally filter out the effect of imaginary components in matrix $\Lambda$ on evolution of $z(t,z_0)$. This equates the degree of stability of equation $\dot{x} = Ax$, which is measured by the maximal real part of eigenvalues of $A$, and equation $\dot{z} = \alpha_1 z$. Furthermore, (2.8) is currently defined explicitly, whereas the definition of its counterpart given in [13] involves numerical evaluation of $w(t)$, $p(w(t))$ and $c(w(t))$.

Albeit that the current version of the auxiliary equation is gaining computational efficiency and widening the scope of the relevant applications, its former counterpart can provide sharper estimates in the common application domain. In fact, our former technique [13] simulates $w(t)$ for equation (1.3) which encapsulates the effects of matrices $A$ and $G_*(t)$ on the time-histories of the corresponding solutions. In contrast, our current technique treats $\|G_*(t)\|$ as a conservative perturbation.

For straight-up further referencing, we write a homogeneous counterpart of (2.8) as follows,

$$\dot{z} = (\alpha_1 + \|G_-(t)\|) z + L(t,z)$$
$$z(t_0, z_0) = z_0 = \|V^{-1}x_0\| \qquad (2.10)$$

Analysis of boundedness/stability of solutions to equations (2.8) or (2.10) is simplified under a condition that brings their linear upper bounds which are defined through utility of the following inequality,

$$L(t,z) \leq l_2(t,\hat{z}) z, \; \forall t \geq t_0, \; \forall z \in [0, \hat{z}], \; \hat{z} > 0 \qquad (2.11)$$

where we assumed that function $l_2 : [t_0, \infty) \times \mathbb{R}_{\geq 0} \to \mathbb{R}_{\geq 0}$ is continuous in $t$ and bounded in $\hat{z}$. For instance, (2.11) can be developed via the application of Lipschitz continuity condition to a scalar nonnegative function $L(t,z)$. Note that definition of $l_2(t,\hat{z})$ is exceedingly simplified in this scalar case.

Furthermore, definition of $l_2(t,\hat{z})$ is also abridged if $L(t,z)$ is a polynomial or piece-wise polynomial function in $z$ which, in former case, can be written as follows, $L(t,z) = \sum_{i=1}^{N} a_i(t) z^i$, where we assume that $a_i(t) \geq 0$. In this case, $l_2(t,\hat{z}) = \sum_{i=1}^{N} a_i(t)(\hat{z})^{i-1}$ since $z^i \leq (\hat{z}^{i-1}) z, \; \forall z \in [0, \hat{z}], \; \hat{z} > 0$.

Clearly, in some cases $l_2$ may merely depend upon $t$ or be a scalar parameter which, however, complies with our further inferences that prompt close-form boundedness/stability criteria and estimation of boundedness/stability regions for equations (1.1) and (1.2).

## 3. Boundedness/Stability of Nonautonomous Nonlinear Systems via Application of Modified Linear Auxiliary Equation

Application of (2.11) to (2.8) yields its linear counterpart that can be written as follows,

$$\dot{Z} = \mu(t, \hat{z})Z + \|F(t)\|, \quad \forall t \geq t_0$$
$$Z(t_0, z_0, \hat{z}) = z_0 = \|V^{-1}x_0\|$$
(3.1)

where $\mu(t, \hat{z}) = \alpha_1 + \|G_-(t)\| + l_2(t, \hat{z})$ and $\|F(t)\| = F_0 \|V^{-1}\eta(t)\|$ with $\|\eta(t)\| = 1$.

To reduce the notation, we adopted in (3.1) and frequently, but always, will set below that
$Z(t, t_0, z_0, \hat{z}) = Z(t, z_0, \hat{z})$.

Clearly, under the conditions of Theorem 1 and inequality (2.11), $\|G_-(t)\|$, $\|F(t)\|$, and $l_2(t, \hat{z})$ are continuous functions which implies that $\mu(t, \hat{z})$ is a continuous in both variables as well.

Thus, under the condition of Theorem 1, (3.1) admits the following solution,

$$Z(t, z_0, \hat{z}) = z_0 Z_h(t, \hat{z}) + F_0 Z_F(t, \hat{z})$$
(3.2)

where function $Z(t, z_0, \hat{z})$ is continuous in $t$ and $z_0$, $Z_h(t, t_0, \hat{z}) = \exp\left(\int_{t_0}^{t} \mu(\tau, \hat{z}) d\tau\right)$,

$\theta(t, \tau, \hat{z}) = \exp\left(\int_{\tau}^{t} \mu(e, \hat{z}) de\right)$ and $Z_F(t, t_0, \hat{z}) = \int_{t_0}^{t} \theta(t, \tau, \hat{z}) \|V^{-1}\eta(\tau)\| d\tau$. As prior, we frequently, but always, are going to use reduce notations, i.e., $Z_h(t, \hat{z}) = Z_h(t, t_0, \hat{z})$ and $Z_F(t, \hat{z}) = Z_F(t, t_0, \hat{z})$. Clearly, $Z_h > 0$ and $\theta > 0$. Thus, the last inequality yields that $Z_F > 0$ as well. We are going to use these inequalities in the proof of Lemma 2.

### 3.1. Stability Criteria and Estimation of Stability Region

To simplify the subsequent references, we will write a homogeneous counterpart to (3.1) as follows,

$$\dot{Z} = \mu(t, \hat{z})Z$$
$$Z(t = t_0, z_0, \hat{z}) = z_0 = \|V^{-1}x_0\|$$
(3.3)

Necessary and sufficient conditions for the stability or asymptotic stability of a linear system directly follow from accessing the behavior of its fundamental solution matrix, see, e.g., a recent paper [42], Lemma 1, and additional references therein. For (3.3), with continuous in $t$ and bounded in $\hat{z}$ function $\mu(t, \hat{z})$, these conditions can be formulated as follows. Equation (3.3) is stable if and only if,

$$\int_{t_0}^{t} \mu(\tau, \hat{z}) d\tau < \infty, \quad \forall \hat{z} \in (0, \hat{z}_b), \forall t_0 \in \mathrm{T}$$
(3.4)

and asymptotically stable if and only if

$$\limsup_{t \to \infty} \left(\int_{t_0}^{t} \mu(\tau, \hat{z}) d\tau\right) = -\infty, \quad \forall \hat{z} \in (0, \hat{z}_B), \forall t_0 \in \mathrm{T}$$
(3.5)

where $\hat{z}_b$ and $\hat{z}_B$ are some values of $\hat{z}$ for which either (3.4) or (3.5) hold.

In turn, due to comparison principal, the stability of a linear scalar equation (3.3) implies stability of the trivial solutions to nonlinear equation (2.10) and, in turn, to (1.2) if $z(t,z_0) \in [0,\hat{z}]$, $\forall t \geq t_0$ since compliance with this condition enables the linearization of (2.10) via application of (2.11). Next lemma provides an upper bound of $Z(t,z_0)$ which, subsequently, aids estimation of $z(t,z_0)$ and, in turn, $\|x(t,x_0)\|$.

**Lemma 1**. Assume that function $\mu(t,\hat{z})$ is continuous in $t$ and bounded in $\hat{z}$, (3.3) is stable for some $\hat{z} \in (0, \hat{z}_{max})$, where $\hat{z}_{max}$ is the maximal value of $\hat{z}$ for which (3.3) is stable that can be infinity, and $z_0 < \hat{z}/Z_s(t_0,\hat{z})$, where $Z_s(t_0,\hat{z}) = \sup_{t \geq t_0} Z_h(t,t_0,\hat{z})$.

Then,
$$\sup_{t \geq t_0} Z(t,z_0,\hat{z}) < \hat{z}, \ \forall t \geq t_0, \ \forall z_0 < \hat{z}/Z_s(t_0,\hat{z}), \ \hat{z} \in (0, \hat{z}_{max}) \qquad (3.6)$$

**Proof.** Clearly, $\sup_{t \geq t_0} Z_h(t,t_0,\hat{z}) = Z_s(t_0,\hat{z}) < \infty$ since $Z_h(t,t_0,\hat{z})$ is continuous in $t$ and (3.3) is stable which yields that $\sup_{t \geq t_0} Z(t,t_0,z_0,\hat{z}) = Z_s(t_0,\hat{z})z_0$ and consequently (3.6) holds □

**Theorem 2.** Assume that assumptions of Theorem 1 are met, $l_2(t,\hat{z})$ is a continuous in $t$ and bounded in $\hat{z}$ function, and equation (3.3) is stable for some $\hat{z} \in (0, \hat{z}_{max})$. Then, the trivial solution to (1.2) is stable and inequality (2.9) takes the form

$$\|x(t,x_0)\| \leq \|V\| z(t,z_0) \leq \|V\| Z(t,z_0,\hat{z}) \leq \|V\| z_0 Z_s(t_0,\hat{z}), \forall t \geq t_0 \in T,$$
$$\forall z_0 < \hat{z}/Z_s(t_0,\hat{z}), \ x_0 \mid \|V^{-1}x_0\| < \hat{z}/Z_s(t_0,\hat{z}), \ \hat{z} \in (0, \hat{z}_{max}) \qquad (3.7)$$

where $Z(t,z_0,\hat{z})$, $z(t,z_0)$ and $x(t,x_0)$ are solutions to equations (3.3), (2.10) and (1.2), respectively, and the region of stability of the trivial solution to (1.2) includes all values of $x_0$ containing into the interior of ellipsoid which is defined as follows,

$$\|V^{-1}x_0\| < \sup_{\hat{z} \in (0, \hat{z}_{max})} (\hat{z}/Z_s(t_0,\hat{z})), \ t_0 \in T \qquad (3.8)$$

**Proof**. Let us show that under conditions of this theorem $z(t,z_0) < \hat{z}$, $\forall t \geq t_0$, $\hat{z} \in (0, \hat{z}_{max})$ if $z_0 < \hat{z}/Z_s$. Pretend in contrary, that $t = t_* > t_0$ is the smallest value of $t$ such that $z(t_*,z_0) = \hat{z}$ under prior conditions. Then, (2.11) and, thus, equation (3.3) holds for $\forall t \in [t_0, t_*]$ which, due to comparison principle [4] and (3.6), implies that $z(t_*,z_0) \leq Z(t_*,z_0,\hat{z}) < \hat{z}$. This contradiction shows that $z(t,z_0) < \hat{z}$, $\forall t \geq t_0$, $\forall z_0 < \hat{z}/Z_s$ which, in turn, enables linearization of (2.10) by (2.11) and brings a linear equation (3.3) that enables our inferences. In order, the application of comparison principal prompts that $z(t,z_0) \leq Z(t,z_0,\hat{z})$, $\forall t \geq t_0$, $z_0 < \hat{z}/Z_s$, $\hat{z} \in (0, \hat{z}_{max})$ which implies that the trivial solution to (2.10) is stable.

The last inequality let to write (2.9) as (3.7) which shows that the stability of (3.3) implies the stability of the trivial solutions to (2.10) and (1.2) under the conditions of this theorem.

Let us now estimate the regions of stability of the trivial solutions to equation (2.10) and (1.2). The solution to (2.10) is stable for all values of $z_0$ which are contingent by the following inequality $z_0 < \sup_{\hat{z} \in (0, \hat{z}_{max})} (\hat{z}/Z_s(t_0,\hat{z}))$.

In sequel, since $z_0 = \|V^{-1}x_0\|$, region of stability of the trivial solution to (1.2) includes all values of $x_0$ meeting condition (3.8) □

Let us recall that (3.8) depends upon $t_0$ - a characteristic property of time – dependent dynamic systems.

**Example 1**. Let us apply inequality (3.8) for the estimation of the stability basin of the planar, nonlinear and nonautonomous equation considered in Section 1.2. Homogeneous counterpart of this planar equation corresponds to nonlinear auxiliary equation (1.22), which, in turn, relates to its linear complement (1.23). For this last equation, we derived in Section 1.2 that $Z_h = \exp\left(\int_{t_0}^{t} s(\tau,\hat{z})d\tau\right) \leq \exp\left(-\overline{s}(\hat{z})(t-t_0)\right)$ which implies that

$$\sup_{t \geq t_0} Z_h(t,t_0,\hat{z}) = Z_s(t_0,\hat{z}) \leq \sup_{t \geq t_0}\left(\exp\left(-\overline{s}(\hat{z})(t-t_0)\right)\right) = 1.$$

Hence, setting in (3.8) that $Z_s(t_0,\hat{z}) = 1$, prompts more conservative but abridged estimation of the stability region for our planar system which can be written as follows, $\|V^{-1}x_0\| < \hat{z}_{\max}$. The last formula was derived also in Section 1.2 by application of a simplified inference.

**Theorem 3**. Assume that assumptions of Theorem 1 are met, $l_2(t,\hat{z})$ is continuous in $t$ and bounded in $\hat{z}$ function, and equation (3.3) is asymptotically stable for some $\hat{z} \in (0,\hat{z}_{\max})$. Then, the trivial solution to (1.2) is asymptotically stable, inequality (3.7) holds where $Z(t,z_0,\hat{z})$, $z(t,z_0)$ and $x(t,x_0)$ are solutions to (3.3), (2.10) and (1.2), respectively, and the region of asymptotic stability of the trivial solution to (1.2) is defined by (3.8).

**Proof**. Literately, under these more conservative conditions, theorem 2 and, hence, (3.7) hold, which affirms the current statement since in this case $\lim_{t \to \infty} Z(t,z_0,\hat{z}) = 0$ and, due to (3.7), $\lim_{t \to \infty} z(t,z_0) = 0$ as well, and, in turn, $\lim_{t \to \infty} x(t,x_0) = 0$ if $z_0, \hat{z}$ and $x_0$ are contingent by conditions of inequality (3.7) and (3.8).

Clearly, in this case (3.8) defines the region of asymptotic stability of the trivial solution to (1.2) □

Next, we present some sufficient stability conditions appending the above statements which can be readily attested in virtually a closed form.

**Corollary 1**. Assume that conditions of Theorem 1 hold, $l_2(t,\hat{z})$ is a continuous in $t$ and bounded in $\hat{z}$ function, and either (3.4) or (3.5) holds as well. Then, inequality (3.7) holds, where $Z(t,z_0,\hat{z})$, $z(t,z_0)$ and $x(t,x_0)$ are solutions to (3.3), (2.10) and (1.2), respectively, and the trivial solution to (1.2) is either stable or asymptotically stable, and the region of stability or asymptotic stability of (1.2) is estimated by (3.8), where in (3.8) $\hat{z}_{\max}$ should be exchanged on either $\hat{z}_b$ or $\hat{z}_B$, respectively.

**Proof**. Really, under the above conditions, (3.3) is either stable or asymptotically stable, which, due to Theorems 2 and 3, assures this statement□

**Corollary 2**. Assume that conditions of Theorem 1 holds, $l_2(t,\hat{z})$ is a continuous in $t$ and bounded in $\hat{z}$ function, and

$$\mu(t,\hat{z}) \leq -\nu(\hat{z}), \forall t \geq t_0, \nu(\hat{z}) > 0, \hat{z} \in (0,\hat{z}_\nu) \tag{3.9}$$

where $\hat{z}_\nu$ is the maximal value of $\hat{z}$ for which the last inequality holds. Then, the trivial solution of (1.2) is asymptotically stable. Inequalities (3.7) and (3.8) hold, where in (3.8) $\hat{z}_{\max}$ should be exchanged on $\hat{z}_\nu$.

**Proof**. In fact, under the above conditions, (3.3) is asymptotically stable. Thus, this statement holds due to Theorem 3□

In sequel, let us set that $\gamma(t_0, \hat{z}) = \limsup_{t \to \infty} (t-t_0)^{-1} \int_{t_0}^{t} m(\tau, \hat{z}) d\tau$, where $m(t, \hat{z}) = \|G_-(t)\| + l_2(t, \hat{z})$

which steers us to

**Corollary 3**. Assume that conditions of Theorem 1 hold, $l_2(t, \hat{z})$ is a continuous in $t$ and bounded in $\hat{z}$ function, and $\phi(t_0, \hat{z}) = \alpha_1 + \gamma(t_0, \hat{z}) < 0$, $\forall \hat{z} \in (0, \hat{z}_\phi)$, $t_0 \in T$, where $\hat{z}_\phi$ is the maximal value of $\hat{z}$ for which last inequality holds. Then, the trivial solution to (1.2) is asymptotically stable, (3.7) and (3.8) hold, where in (3.8) $\hat{z}_{\max}$ should be exchanged on $\hat{z}_\phi$.

**Proof**. In fact, condition, $\phi(t_0, \hat{z}) < 0$ yields that the Lyapunov exponent of solutions to (3.3) is negative, which assures that (3.3) is asymptotically stable [1]. Hence, due to theorem 3, (3.7) and (3.8) hold, where $Z(t, z_0, \hat{z})$, $z(t, z_0)$ and $x(t, x_0)$ are solutions to (3.3) and (2.10) and (1.2), respectively, which attests this statement □

**Remark**. Let us represent the last stability criterion in an abridged form. Pretend that

$\|G_-(t)\| = g_{av}(t_0) + \tilde{g}(t_0, t)$, where $g_{av} = \limsup_{t \to \infty} (t-t_0)^{-1} \int_{t_0}^{t} \|G_-(\tau)\|(\tau) d\tau$, $\tilde{g} = \|G_-(t)\| - g_{av}$. Also

suppose that $l_2(t, \hat{z}) = l_{2av}(t_0, \hat{z}) + \tilde{l}_2(t_0, t, \hat{z})$, $l_{2av}(t_0, \hat{z}) = \limsup_{t \to \infty} (t-t_0)^{-1} \int_{t_0}^{t} l_2(t_0, \tau, \hat{z}) d\tau$,

$\tilde{l}_2 = l_2(t, \hat{z}) - l_{2av}(t_0, \hat{z})$. Then, it is follows from definition of $\gamma$ that $\gamma(t_0, \hat{z}) = g_{av}(t_0) + l_{2av}(t_0, z_0)$ which implies that stability condition of Lemma 3 can be written as follows,

$$\alpha_1 + g_{av}(t_0) + l_{2av}(t_0, \hat{z}) < 0, t_0 \in T$$

Note that the above statements grant stability criteria of the trivial solutions to (1.2) through the application of readily verifiable and frequently closed form formulas if $f(t, x)$ is a piecewise polynomial in $x$ since in this case functions $L(t, z)$ and $l_2(t, \hat{z})$ can be defined in analytical form. In contrast, the utilities of either (1.6) or (1.8) require simulations of $w(t)$ and subsequent estimation of the corresponding parameters or function from the simulated data. The solution to such inverse problems can be sensitive to the variation of some of their parameters, which, for instance, can include the length of time-interval, etc. Consequently, to our knowledge, these problems have not been attended yet for practically feasible systems.

### 3.2. Boundedness Criteria and Estimation of Region of Boundedness

Next, let us develop some criteria of boundedness of solutions to equation (1.1) by assessing the relevant behavior of solutions to equation (3.1). These criteria frequently resemble the criteria of stability/asymptotic stability of forced solutions to nonhomogeneous systems. To our knowledge, such criteria and especially the estimation of regions of boundedness of solutions to multidimensional, time – varying nonlinear systems are rare in the current literature. Yet, there is some relation between our boundedness criteria and the concept of input-to-state stability developed by E. D. Sontag under more restrictive conditions, see [43], [44] and [4].

A solution to (3.1), $Z(t, t_0, z_0, \hat{z}) < \infty$, $\forall t \geq t_0$, $\hat{z} \in (0, \hat{z}_1)$ if (3.3) is either stable or asymptotically stable and

$$\rho(t_0, \hat{z}) = \sup_{t \geq t_0} Z_F(t, t_0, \hat{z}) < \infty, \hat{z} \in (0, \hat{z}_\rho), \forall t_0 \in T \tag{3.10}$$

where $\hat{z}_\rho$ is the maximal value of $\hat{z}$ for which (3.10) holds, and $\hat{z}_1 = \min(\hat{z}_{\max}, \hat{z}_\rho)$.

Firstly, let us present an analog to Lemma 1 for the nonhomogeneous equation (3.1) as follows,

**Lemma 2**. Assume that function $\mu(t,\hat{z})$ is a continuous in $t$ and bounded in $\hat{z}$, (3.3) is stable and

$$z_0 < (\hat{z} - F_0\rho(t_0,\hat{z}))/Z_s(t_0,\hat{z}) = z_F(t_0,\hat{z}), \ \hat{z} \in (0,\hat{z}_1) \tag{3.11}$$

Then,

$$\sup_{t \geq t_0} Z(t,z_0,\hat{z}) < Z_s(t_0,\hat{z})z_0 + F_0\rho(t_0,\hat{z}) < \hat{z}, \ \forall t \geq t_0, \ \forall z_0 < z_F \tag{3.12}$$

where, as prior, $Z_s(t_0,\hat{z}) = \sup_{t \geq t_0} Z_h(t,t_0,\hat{z})$.

**Proof.** As is mentioned in Lemma 1, $\sup_{t \geq t_0} Z_h(t,t_0,\hat{z}) = Z_s(t_0,\hat{z}) < \infty$ under conditions of this statement. Then,

$\sup_{t \geq t_0} Z(t,t_0,z_0,\hat{z}) \leq z_0 \sup_{t \geq t_0} Z_h(t,t_0,\hat{z}) + \sup_{t \geq t_0} Z_F(t,t_0,\hat{z}) = z_0 Z_s(t_0,\hat{z}) + F_0\rho(t_0,\hat{z})$ since $Z_h(t_0,t,\hat{z}) \geq 0$

and $Z_F(t_0,t,\hat{z}) \geq 0$ which yields (3.11)□

The next statement assesses the boundedness of solutions to the nonhomogeneous equations (2.8) and, in turn, (1.1).

**Theorem 4.** Assume that the assumptions of Theorem 1 hold, equation (3.3) is stable for some $\hat{z} \in (0,\hat{z}_{\max})$, inequalities (3.10) and (3.11) hold, and $l_2(t,\hat{z})$ is a continuous in $t$ and bounded in $\hat{z}$ function. Then,

(I) $\|x(t,x_0)\| \leq \infty, \forall t \geq t_0$ if these solutions are emanated from the interior of the ellipsoid that is defined as follows,

$$x_0 \mid \|V^{-1}x_0\| < \sup_{\hat{z} \in (0,\hat{z}_1)} z_F(t_0,\hat{z}) \tag{3.13}$$

(II) inequality (2.9) takes the form,

$$\|x(t,x_0)\| \leq \|V\| \|z(t,z_0)\| \leq \|V\| Z(t,z_0,\hat{z}) = \|V\|(z_0 Z_h(t,z_0) + F_0 Z_F(t,z_0)) \leq$$
$$\|V\|(z_0 Z_s(t_0,\hat{z}) + F_0\rho(t_0,\hat{z})); \forall t \geq t_0, \ x_0 \mid \|V^{-1}x_0\| < z_F, \ z_0 < z_F, \ \hat{z} \in (0,\hat{z}_1), \ \forall t_0 \in \mathrm{T} \tag{3.14}$$

where $Z(t,t_0,z_0,\hat{z})$, $z(t,z_0)$ and $x(t,t_0,x_0)$ are solutions to (3.1) and (2.8) and (1.1), respectively, and

(III) if (3.3) is asymptotically stable and other conditions of this statement hold, then,

$$\limsup_{t \to \infty} \|x(t,x_0)\| \leq F_0 \|V\| \rho(t_0,\hat{z}), \ x_0 \mid \|V^{-1}x_0\| < z_F, \hat{z} \in (0,\hat{z}_1), \ \forall t_0 \in \mathrm{T} \tag{3.15}$$

**Proof**. Firstly, we show that under the conditions of this theorem $z(t,z_0) < \hat{z}, \ \forall t \geq t_0, \ \hat{z} \in (0,\hat{z}_1)$ if $z_0 < z_F$. In fact, pretend, in contrary, that $t = t_* > t_0$ is the smallest value of $t$ such that $z(t_*,z_0) = \hat{z}$ under prior conditions. Then, (2.11) and, thus, equation (3.1) holds for $\forall t \in [t_0,t_*]$ which, due to comparison principle and (3.11), implies that $z(t_*,z_0) \leq Z(t_*,z_0,\hat{z}) < \hat{z}$. This contradiction prompts that

$z(t,z_0) < \hat{z}, \ \forall t \geq t_0, \ \forall z_0 < z_F, \hat{z} \in (0,\hat{z}_1)$ which, in turn, enables a linearization of (2.8) with the utility of (2.11) that, subsequently, yields equation (3.1). In order, the application of comparison principal returns that $z(t,z_0) \leq Z(t,z_0,\hat{z}) < \infty, \ \forall t \geq t_0, \ z_0 < z_F, \ \hat{z} \in (0,\hat{z}_1)$ since (3.3) is stable and (3.11) holds.

Clearly, the last inequality lets us write (2.9) as (3.14), which demonstrates that the solutions of (1.1) are bounded within the region specified by (3.13).

Lastly, (3.15) directly follows from (3.14) if (3.3) is asymptotically stable since in this case $\lim_{t\to\infty} Z_h(t, x_0) = 0$ □

Next, two corollaries present simplified sufficient conditions under which the presumptions of Theorem 4 are met.

**Corollary 5**. Assume that the conditions of Theorem 1 and inequalities (3.9), (3.10) and (3.11) hold, and $l_2(t, \hat{z})$ is a continuous in $t$ and bounded in $\hat{z}$ function. Then, the solutions to (1.1) are bounded in the region that is defined by (3.13), inequality (3.14) holds, and
$$\limsup_{t\to\infty} \|x(t, x_0)\| \leq F_0 \|V\| \|V^{-1}\| / v, \ \forall x_0 \mid \|V^{-1} x_0\| = z_F$$

Note that this statement assumes that in (3.13) and (3.14) $\hat{z}_1 = \min(\hat{z}_v, \hat{z}_\rho)$.

**Proof**. In fact, the condition of this statement implies the asymptotic stability of (3.3) for $\forall \hat{z} \in (0, \hat{z}_v]$ which, together with other assumptions, yields that the conditions of Theorem 4 hold. In turn, this implies that (3.14) holds which yields that $\limsup_{t\to\infty} \|x(t, x_0)\| \leq F_0 \|V\| \sup_{t \geq t_0} Z_F(t, x_0)$. Next,

$$Z_F(t, \hat{z}) = \int_{t_0}^{t} \theta(t, \tau, \hat{z}) \|V^{-1} \eta(\tau)\| d\tau \leq \int_{t_0}^{t} \|\theta(t, \tau, \hat{z})\| \|V^{-1}\| \|\eta(\tau)\| d\tau \leq \|V^{-1}\| \int_{t_0}^{t} e^{-v(t-\tau)} d\tau \leq \|V^{-1}\|/v \ \ □$$

In the sequel, let us assume that the Lyapunov exponent of solutions to (3.3), i.e., $\phi(\hat{z}) = \alpha_1 + \gamma(\hat{z}) < 0$, $\hat{z} \in (0, \hat{z}_\phi)$, $\hat{z}_\phi > 0$, see corollary 3 for additional details. This implies that $Z_h(t, z_0, \hat{z}) \leq D(\vartheta, \hat{z}) e^{-\chi(t-t_0)}$ and, due to (1.15), that $\theta(t, \tau, \hat{z}) \leq D(\vartheta, \hat{z}) e^{-\chi(t-\tau)}$, where, $-\chi = \phi + \vartheta$, $\chi > 0$, $\vartheta > 0$ is a small value and constants $D(\vartheta, \hat{z}) > 0$ [6]. This leads to

**Corollary 6**. Assume that the conditions of Theorem 1, inequalities (3.10) and (3.11) hold, $\phi(t_0, \hat{z}) < 0$, $\forall \hat{z} \in (0, \hat{z}_\phi)$, and function $l_2(t, \hat{z})$ is a continuous in $t$ and bounded in $\hat{z}$. Then, (3.14) and (3.15) hold and $\lim_{t\to\infty} \|x(t, x_0)\| \leq F_0 \|V\| \|V^{-1}\| D(\vartheta, \hat{z})/\chi, \ \forall x_0 \mid \|V^{-1} x_0\| = z_F$.

Note that in this statement assumes that in (3.13) and (3.14) $\hat{z}_1 = \min(\hat{z}_\phi, \hat{z}_\rho)$.

**Proof**. In fact, in this case (3.3) is asymptotically stable [1] and the conditions of theorem 4 hold. Thus, (3.14) yields that $\limsup_{t\to\infty} \|x(t, x_0)\| \leq F_0 \|V\| \sup_{t \geq t_0} Z_F(t, x_0)$. In turn,

$$Z_F(t, t_0, \hat{z}) = \int_{t_0}^{t} \theta(t, \tau, \hat{z}) \|V^{-1} \eta(\tau)\| d\tau \leq \int_{t_0}^{t} \|\theta(t, \tau, \hat{z})\| \|V^{-1}\| \|\eta(\tau)\| d\tau \leq$$

$$\|V^{-1}\| D(\vartheta, \hat{z}) \int_{t_0}^{t} e^{-\chi(t-\tau)} d\tau \leq \|V^{-1}\| D(\vartheta, \hat{z})/\chi \ □$$

**Example 2**. To abridge our derivation, we apply this methodology to estimation of the region of boundedness of solutions to (1.21), which, in turn, prompts estimation of the corresponding region for the planar system considered in Section 1.2. Yet, (1.21) is a scalar but nonintegrable equation and direct estimation of its boundedness region requires repeated simulations of this equation. Thus, to simplify our further inferences, we bring an autonomous counterpart to (1.21) which can be written as follows

$$\begin{aligned} \dot{z}_1 &= (-\alpha_* + G_s) z_1 + \sigma z_1^3 + \tilde{F} \\ z_1(t_0, z_0) &= z_0 = \|V^{-1} x_0\| \end{aligned} \tag{3.16}$$

where, as prior, we assume that $G_s = \sup_{\forall t \geq t_0} \|G_-(t)\| < \infty$ and $\|F(t)\| \leq F_0$, $\|V^{-1}\| = \tilde{F} = const$. Thus, due to comparison principle, $z(t, z_0) \leq z_1(t, z_0)$, $\forall t \geq t_0$, where $z(t, z_0)$ is a solution to (1.21). In turn, the boundedness of solutions to the autonomous and integrable equation (3.16) is determined by the location of its readily accessible fixed solutions.

As in Section 1.2, to simplify the assessment of the behavior of solutions to (3.16), we consider its linear counterpart

$$\dot{Z}_1 = -\overline{s}(\hat{z})Z_1 + \tilde{F}, \; \forall t \geq t_0 \qquad (3.17)$$
$$Z_1(t_0, z_0, \hat{z}) = z_0 = \|V^{-1}x_0\|$$

where $-\overline{s}(\hat{z}) = -\alpha_* + G_s + (\hat{z})^2 \sigma$. Thus, (3.17) can be viewed as a nonhomogeneous analog to (1.25).

Furthermore, let us write a nonhomogeneous counterpart to (1.23) as follows

$$\dot{Z}_2 = s(t, \hat{z})Z_2 + \tilde{F}, \; \forall t \geq t_0, \; z_0 < \hat{z} \qquad (3.18)$$
$$Z_2(t_0, z_0, \hat{z}) = z_0 = \|V^{-1}x_0\|$$

Thus, due to comparison principle, $z_1(t, z_0) \leq Z_1(t, z_0, \hat{z})$, $\forall t \geq t_0$ and $Z_2(t, z_0, \hat{z}) \leq Z_1(t, z_0, \hat{z})$, $\forall t \geq t_0$.

Clearly, the general solutions to the linear scalar equations (3.17) and (3.18) take the form, $Z_i(t, z_0, \hat{z}) = z_0 Z_{i,h}(t, z_0, \hat{z}) + \tilde{F} Z_{i,F}(t, z_0, \hat{z})$, $i = 1, 2$, $\forall t \geq t_0$, where $Z_{i,h}(t, z_0, \hat{z})$, $i = 1, 2$ are solutions to (1.25) and (1.23) with $z_0 = 1$ and $Z_{i,F}(t, z_0, \hat{z})$ are components of forced solution to the corresponding equations. Let us recall that $\sup_{\forall t \geq t_0} Z_{1,h}(t, z_0, \hat{z}) = Z_{1,s}(z_0, \hat{z}) = \sup_{\forall t \geq t_0} \exp(-\overline{s}(\hat{z})(t - t_0)) = 1$, see Section 1.2.

To estimate the right – side of (3.11), we recall that $Z_s(z_0, \hat{z}) = \sup_{\forall t \geq t_0} Z_{2,h}(t, z_0, \hat{z}) \leq \sup_{\forall t \geq t_0} Z_{1,h}(t, z_0, \hat{z}) = 1$.

In sequence, $\theta(t, \tau, \hat{z}) = \exp\left(\int_\tau^t s(e, \hat{z})de\right) \leq \exp\left(\int_\tau^t -\overline{s}(\hat{z})de\right) \leq \exp(-\overline{s}(\hat{z})(t-\tau))$ which, with use of (3.10), let us bound $\rho(t_0, \hat{z})$ as follows

$$\rho(t_0, \hat{z}) = \sup_{\forall t \geq t_0} \int_{t_0}^t \theta(t, \tau, \hat{z}) \|V^{-1}\eta(\tau)\| d\tau \leq \|V^{-1}\| \sup_{\forall t \geq t_0} \int_{t_0}^t \theta(t, \tau, \hat{z}) d\tau \leq$$
$$\|V^{-1}\| \sup_{\forall t \geq t_0} \int_{t_0}^t \exp(-\overline{s}(\hat{z})(\tau - t_0)) d\tau \leq (\|V^{-1}\|/s(\hat{z})) \sup_{\forall t \geq t_0} (1 - \exp(-s(\hat{z})(t - t_0))) = \|V^{-1}\|/\overline{s}(\hat{z}) \qquad (3.19)$$

Then, using the developed inequalities, we bound the right – side of (3.13) as follows,

$$\sup_{\hat{z} \in (0, \hat{z}_1)} z_F(\hat{z}) = \sup_{\hat{z} \in (0, \hat{z}_1)} \left(\left(\hat{z} - F_0 \rho(t_0, \hat{z})\right)/Z_s(t_0, \hat{z})\right) \leq \sup_{\hat{z} \in (0, \hat{z}_1)} \left(\hat{z} - \tilde{F}/\overline{s}(\hat{z})\right)$$

where in this case $\hat{z}_1 = \hat{z}_*$ and $\hat{z}_*$ is the smallest positive root of equation $\overline{s}(\hat{z}) = 0$.

Next, we inference that $\sup_{\hat{z}\in(0,\hat{z}_1)} \left(\hat{z} - \tilde{F}/\overline{s}(\hat{z})\right)$ is reached for $\hat{z} = \hat{z}_s < z_*$ since $\lim_{\hat{z}\to z_*} 1/\overline{s}(\hat{z}) = \infty$. Note that $z_s$ can be readily found numerically. Hence, the region of boundedness of solutions to equation (3.17) and, thus, to (3.16) and, subsequently, (1.21) can be estimated as follows,

$$z_0 < \hat{z}_s - \tilde{F}/\overline{s}(\hat{z}_s) = R$$

At this point, we consider a numerical example with $\alpha_* = 0.3$ and $G_s = \sigma = \tilde{F} = 0.1$. For this set of parameters, (3.16) has three real fixed solutions in which two solutions are positive and equal $0.618$ and one. Furthermore, in this case, the smallest positive fixed solution is stable whereas the largest - is unstable. Thus, the solutions to (3.16) are bounded and approach the stable fixed solution as $t \to \infty$ if $z_0 \in [0,1)$, and approach infinity as $t \to \infty$ if $z_0 > 1$. Yet, for this set of parameters, $z_* = \sqrt{2}$, $z_s = 0.848$ and $R = 0.07$. The last value provides a rather conservative estimate of the boundedness interval of solutions to (3.16) since our current derivation rest on a sequence of upper bounds which simplifies our inference. More accurate estimates can be obtained through direct simulation of (3.13). Nonetheless, the above reasoning authenticates our theory.

## 4. Modified Nonlinear Auxiliary Equation and Boundedness/ Stability of Nonautonomous and Nonlinear Systems

Analysis of solutions to the scalar nonlinear auxiliary equations (2.8) and (2.10) conveys the less conservative criteria of boundedness/stability of the initial nonlinear and multidimensional systems and naturally enfolds the estimation of their trapping/stability regions.

In general, (2.8) is a nonintegrable, scalar and nonlinear equation with time-varying coefficients and the behavior of its solutions can be readily assessed either in simulations or by laying out some integrable nonlinear and autonomous equations bounding (2.8) from above and below. The former approach is substantially abridged due to the following

**Theorem 5**. Assume that equation (2.8) possess a unique solution $z(t, z_0)$, $\forall t \geq t_0$, $z_0 \geq 0$ and $z(t, z_0')$ and $z(t, z_0'')$ are solutions to (2.8) with $0 \leq z_0' \leq z_0''$. Then, $z(t, z_0') \leq z(t, z_0'')$, $\forall t \geq t_0$.

**Proof**. In fact, the solutions to (2.8) do not intersect in $t \times z$-plane due to uniqueness property of this equation □

Application of the above statement grants some boundedness/stability criteria for equations (1.1) and (1.2), respectively. Firstly, let us define a set of centered at zero concentric ellipsoids, $E(z) \subset \mathbb{R}^n$ as follows,

$$E(z): \|V^{-1}x\| = z \geq 0 \qquad (4.1)$$

Also, we assume that $\partial E(z) \subset \mathbb{R}^{n-1}$ defines the boundaries of these ellipsoids and $E_-(z) = E(z) - \partial E(z)$. This prompts the following,

**Theorem 6**. Assume that equations (1.2) and (2.10) possess unique solutions for $\forall x_0 \in H \subset \mathbb{R}^n$ and $\forall z_0 = \|V^{-1}x_0\|$, and that the trivial solution to (2.10) is either stable or asymptotically stable. Then the trivial solution to (1.2) is stable/asymptotically stable as well, respectively. Furthermore, assume that the interval $[0, \overline{z})$ defines the stability region of the stable/asymptotically stable trivial solution to (2.10). Then set $E_-(\overline{z})$ is contained in stability region of (1.2).

**Proof**. Really, the proof of this statement follows from inequality (2.9), where it is presumed that $x(t, x_0)$ and $z(t, z_0)$ are solutions to equations (1.2) and (2.10), respectively □

**Theorem 7**. Assume that equations (1.2) and (2.8) possess unique solutions for $\forall x_0 \in H \subset \mathbb{R}^n$ and $\forall z_0 = \|V^{-1}x_0\|$, and interval $[0, \overline{z}]$ defines a trapping region of (2.8) about $z \equiv 0$, i.e., $\|z(t, z_0)\| \leq \overline{z}$,

$\forall z_0 \in [0, \overline{z}]$, $\forall t \geq t_0$. Then, the ellipsoid $E(\overline{z})$ is included into the trapping region of equation (1.1) about zero, i.e., $\|x(t, x_0)\| < \infty$, $\forall x_0 \in E(\overline{z})$, $\forall t \geq t_0$.

**Proof**. As formerly, the proof of this statement follows from inequality (2.9), where $x(t, x_0)$ and $z(t, z_0)$ are solutions to equations (1.1) and (2.8), respectively □

Thus, the problem of estimating the trapping/stability regions about point $x \equiv 0$ of multidimensional equations (1.1) or (1.2) is reduced to reckoning the threshold value $\overline{z}$ splitting the interval of initial values of solutions to (2.8) or (2.10) in two parts associated with qualitatively distinct behavior of such solutions on long time-intervals. Then the boundary of the trapping/stability region for either (1.1) or (1.2) is given by the following formula, $\overline{z} = \|V^{-1} x_0\|$. In turn, the task of simulating the threshold value is markedly simplified since $z(t, z_0)$ is monotonically increases in $z_0$, $\forall t \geq t_0$ due to theorem 5.

Yet, the structure of solutions to (2.8) or (2.10) and, in turn, to their multidimensional counterparts (1.1) and (1.2) can be further divulged through the analysis of solutions to some simplified equations which bound (2.8) or (2.10). In this context, we acknowledge again that under (2.11) the solutions to linear equations (3.1) and (3.3) bounds from above the corresponding solutions to equations (2.8) and (2.10) that are stemmed from the same initial values.

A distinct stability criterion can be recouped using the Massera, Chetaev, and Malkin theorem, which replaces (2.11) by a complementary condition, see [1] for contemporary references and citations. In this regard, we assume that $L(t, z) = L_1(t) z + L_r(t, z)$, where $L_1 : [t_0, \infty) \to \mathbb{R}_{\geq 0}$ is continuous in $t$ and $L_r : [t_0, \infty) \times \mathbb{R}_{\geq 0} \to \mathbb{R}_{\geq 0}$ is continuous in both variables and also Lipschitz in $z$, and

$$L_r(t, z) \leq l_3 z^r, \quad \forall t \geq t_0, \; r > 1, \; \forall z \in [0, \tilde{z}], \; l_3, \tilde{z} > 0 \tag{4.2}$$

Then (2.10) can be written as

$$\dot{z} = \psi(t) z + L_r(t, z)$$
$$z(t_0) = z_0 = \|V^{-1} x_0\| \tag{4.3}$$

where $\psi(t) = \alpha_1 + \|G_-(t)\| + L_1(t)$.

Next, let $\varsigma(t, \tau) = \int_\tau^t \psi(s) ds$, where

$$\varsigma(t, \tau) \leq \varsigma_0 - \varsigma_1(t - \tau) + \varsigma_2 \tau, \; \forall t \geq \tau \geq 0, \; \varsigma_0, \varsigma_1 > 0, \; \varsigma_2 \geq 0, \tag{4.4}$$

and

$$(r-1)\varsigma_1 \geq \varsigma_2 \tag{4.5}$$

This comprises the following,

**Theorem 8**. Assume that equations (1.2) and (2.10) possess unique solutions for $\forall x_0 \in H \subset \mathbb{R}^n$ and $\forall z_0 \geq 0$, and that (4.2) is contended with sufficiently small $l_3$ and (4.4) and (4.5) hold as well. Then, the trivial solution to (4.3) is asymptotically stable, which, in turn, due to the comparison principle, implies that the trivial solution to the corresponding equation (1.2) is asymptotically stable as well.

**Proof**. Indeed, the proof of this statement directly follows from application of the mentioned above theorem to equation (4.3) which represents equation (2.10) in this case. Then, this statement is assured due to application of the comparison principle □

Some additional approaches to analysis of the structure of solutions to a scalar auxiliary equation were outlined in [13]. Below we complement and apply some of these techniques to a modified auxiliary equation (2.8) which is developed in this paper. To this intend, let us write two scalar and autonomous equations as follows,

$$\dot{z}_2 = \kappa_+ z_2 + L_+(z_2) + F_0$$
$$z_2(t_0) = z_0 = \|V^{-1} x_0\| \tag{4.6}$$

$$\dot{z}_3 = \kappa_- z_3 + L_-(z_3) + F_-$$
$$z_3(t_0) = z_0 = \|V^{-1} x_0\|$$
(4.7)

where $\kappa_+ = \alpha_1 + \sup_{t \geq t_0} \|G_-(t)\|$, $\kappa_- = \alpha_1 + \inf_{t \geq t_0} \|G_-(t)\|$, $L_+(z) = \sup_{t \geq t_0} L(t, z)$, $L_-(z) = \inf_{t \geq t_0} L(t, z)$, $z \geq 0$, $F_- = F_0 \inf_{t \geq t_0} \|\eta(t)\|$.

In turn, we assume in the reminder of this section that equations (4.6) and (4.7) possess unique solutions, $z_i(t, z_0)$, $\forall t \geq t_0$, $\forall z_0 \geq 0$, $i = 2, 3$.

Clearly, the right sides of the last two scalar, autonomous and integrable equations bound from above and below the right side in (2.8), which, due to the comparison principle, implies that, $z_3(t, z_0) \leq z(t, z_0) \leq z_2(t, z_0)$, $\forall t \geq t_0$, where $z(t, z_0)$ is a solution to (2.8). In turn, the structures of solutions to the autonomous and scalar equations (4.6) and (4.7) are determined by the location and stability of their fixed solutions. Application of such reasoning to (4.6) yields sufficient conditions for the boundedness and stability of solutions to (2.8) and (2.10), respectively, which, in turn, embraces the corresponding statements for solutions to (1.1) and (1.2). Conversely, resolving the behavior of solutions to (4.7), brings the necessary conditions for the boundedness and stability of solutions to (2.8) and (2.10), respectively.

First, we presume below that $\kappa_+ < 0$ and $\kappa_- < 0$ since, otherwise, the right sides of (4.6) and (4.7) become positive which entails that $\lim_{t \to \infty} z_i(t, z_0) = \infty$, $\forall z_0 \geq 0$ $i = 2, 3$. The last conditions imply that in these cases our bounds, i.e., $z_i(t, z_0)$ become over conservative if values of $t$ become relatively large.

Next, we recap the akin sufficient conditions in the following,

**Theorem 9**. Presume that equations (1.2), (2.8), (2.10), and (4.6) possess unique solutions for $\forall x_0 \in H \subset \mathbb{R}^n$ and $\forall z_0 = \|V^{-1} x_0\|$. Additionally, assume (I) that $F_0 = 0$, $\kappa_+ < 0$ and $z_c \neq 0$ is a unique unstable fix solution to (4.6). Then the trivial solution to (2.10) is asymptotically stable, i.e., $\lim_{t \to \infty} z(t, z_0) = 0$, $\forall z_0 \in [0, z_c)$ and, consequently, $z_c$ bounds from above the stability region of the trivial solution to (2.10). This, in turn, implies that the trivial solution to (1.2) is asymptotically stable and $\lim_{t \to \infty} \|x(t, x_0)\| = 0$, $\forall x_0 \in E_-(z_c)$, i.e. $E_-(z_c)$ is enclosed into stability region of the trivial solution to (1.2).

Next, let us assume that $F_0 > 0$, $\kappa_+ < 0$ and (II) (4.6) has two fix solutions $0 < z_{c2} < z_{c1}$ corresponding to simple roots of its right side. Then, $z_{c1}$ and $z_{c2}$ are unstable and stable fix solutions to (4.6), respectively. Furthermore, a solution to (2.8), $z(t, z_0) \leq z_{ci}$, $\forall z_0 \leq z_{ci}$, $i = 1, 2$ and $\lim_{t \to \infty} z_2(t, z_0) = z_{c2}$, $\forall z_0 < z_{c1}$. This, consequently, yields that

$$\|x(t, x_0)\| \leq \begin{cases} \|V\| z_{c1}, & \forall x_0 \in E(z_{c1}) \\ \|V\| z_{c2}, & \forall x_0 \in E(z_{c2}) \end{cases}, \forall t \geq t_0,$$

$$\limsup_{t \to \infty} \|x(t, x_0)\| \leq \|V\| z_{c2}, \forall x_0 \in E_-(z_{c1})$$

Lastly, let us assume that $F_0 > 0$, $\kappa_+ < 0$ and (III) $z_{c1} = z_{c2} = z_C$. Then, a solution to (2.8), $z(t, z_0) \leq z_C$, $\forall z_0 \leq z_C$ which, in turn, implies that $\|x(t, x_0)\| \leq \|V\| z_C$, $\forall x_0 \in E_-(z_C)$, $\forall t \geq t_0$.

**Proof**. The proof of this statement immediately follows from the qualitative analysis of the behavior of solutions to (4.6) which bound from above solutions to either (2.8) or (2.10) and the subsequent application of inequality (2.9) that, in turn, bounds time-histories of $\|x(t, x_0)\|$.

Let us, for instance, show that $z_{c1}$ and $z_{c2}$ are unstable and stable fixed solutions to (4.6). In fact, we get from (4.6) that $\dot{z}_{2|z_2=0} = F_0 > 0$ which implies that continuous right side of (4.6) is a positive function for $\forall z < z_{c2}$ since $z_{c2}$ is a simple and smallest root of equation, $\kappa_+ z_2 + L_+(z_2) + F_0 = 0$. Hence, $z_{c2}$ is a stable and, in turn, $z_{c1}$ is an unstable fix solution to (4.6)□

Application of similar reasoning to equation (4.7) leads to necessary conditions for the boundedness and stability of solutions to (2.8) and (2.10), respectively, which can aid simulations of these equations. Nonetheless, these conditions do not directly endorse the consistent properties of solutions to (1.1) or (1.2). Yet, the utility of the lower bounds of solutions to (2.8) can simplify numerical resembling of the threshold initial value of solutions to this equation. In fact, assume that $\hat{z}_0$ is a threshold initial value for (4.7) such that, for instance, $\lim_{t \to \infty} z_3(t, z_0) < \infty$, $\forall z_0 \leq \hat{z}_0$ and $\lim_{t \to \infty} z_3(t, z_0) = \infty$, $\forall z_0 > \hat{z}_0$. Then, $\lim_{t \to \infty} z(t, z_0) = \infty$, $\forall z_0 > \hat{z}_0$ as well. Thus, $\hat{z}_0$ yields an upper estimate of the actual threshold initial value of solutions to (2.8) which can be accessed through analytical reasoning and aid the approximate resembling of the boundaries of trapping/stability regions for this equation.

Next, we endorse an additional lower bound for solutions to (2.10) through the utility of the integrable Bernoulli equation which is bounded from below the right side of (2.10). In fact, assume that function $L(t, z)$ in this equation can be written as follows, $L(t, z) = L_1(t) z + L_n(t) z^n + P(t, z)$, where continuous, $L_1(t), L_n(t) > 0, \forall t \geq t_0, 1 < n \in \mathbb{N}$, and function $P(t, z) \geq 0$, $z \geq 0, \forall t \geq t_0$ is continuous and Lipschitz in $z$. This enfolds the following Bernoulli equation which assumes analytical solutions,

$$\dot{z}_4 = \left(\alpha_1 + \|G_-(t)\| + L_1(t)\right) z_4 + L_2(t) z_4^n$$
$$z_4(t_0, z_0) = z_0 = \|V^{-1} x_0\|$$
(4.8)

Clearly, $z_4(t, z_0) \leq z(t, z_0), \forall t \geq t_0$, and $z_4(t, z_0)$ should be close to $z(t, z_0)$ if $P(t, z)$ is sufficiently small. Thus, the threshold initial value procured in the testing of solutions to the integrable equation (4.8) bounds from above actual threshold value which can be retrieved in simulations of (2.10).

## 5. Simulations

This section applies the developed methodology to decipher the evolutions of solutions' norms and estimate the regions of stability/ boundedness of two nonlinear systems with time – dependent nonperiodic coefficients which are common in various applications [40]. These systems comprise of two coupled Van der – Pol – like or Duffing - like oscillators with variable coefficients. In both cases, the direct utility of our prior technique [13] is compromised since in these cases, $\lim_{t \to \infty} c(t) = \infty$.

### 5.1 Coupled System of Van der – Pol – like Oscillators

Equation (2.1) can be turned into such a system if, for instance, we assume that

$$A = \begin{pmatrix} 0 & 1 & 0 & 0 \\ -(\omega_1^2 + d) & -\alpha_1 & d & 0 \\ 0 & 0 & 0 & 1 \\ d & 0 & -(\omega_2^2 + d) & -\alpha_2 \end{pmatrix}, \quad G_* = \begin{pmatrix} 0 & 0 & 0 & 0 \\ -g_{21}(t) & 0 & 0 & 0 \\ 0 & 0 & 0 & 0 \\ 0 & 0 & -g_{43}(t) & 0 \end{pmatrix},$$
(5.1)
$$F_*(t) = \begin{pmatrix} 0 & F_{*1}(t) & 0 & F_{*2}(t) \end{pmatrix}^T$$

$$f_* = \begin{pmatrix} 0 & -\mu_1 x_2^3 & 0 & -\mu_2 x_4^3 \end{pmatrix}^T$$
(5.2)

where $\alpha_1 = 0.4$, $\alpha_2 = 0.2$, $\omega_1^2 = 1$, $\omega_2^2 = 4$, $g_{21} = a_1 \sin r_1 t + a_2 \sin r_2 t$, $g_{43} = b_1 \sin s_1 t + b_2 \sin s_2 t$, $a_1 = a_2 = b_1 = b_2 = 0.1$, $r_1 = 3.14$, $r_2 = 6.15$, $s_1 = 3.1$, $s_2 = 6.28$, $F_{*i} = F_{0i} \sin q_i t$, $i = 1, 2$, $q_1 = 5.43$, $q_2 = 10$. Note that most but all our simulations assume that $\mu_1 = \mu_2 = 18.2$.

Next, let $V = [v_{ij}]$, $i, j = 1, \ldots, 4$ is the eigenvector matrix for matrix $A$ that is defined above and in equation (2.2) $\Lambda = diag(\lambda_1, \lambda_1^*, \lambda_2, \lambda_2^*)$, $\lambda_k = \alpha_k + i\beta_k$, $\lambda_k^* = \alpha_k - i\beta_k$, $k = 1, 2$, $\alpha_1 \geq \alpha_2$,

$$f = V^{-1} \left( 0 \quad -\mu_1 \left( \sum_{k=1}^{4} v_{2k} y_k \right)^3 \quad 0 \quad -\mu_2 \left( \sum_{k=1}^{4} v_{4k} y_k \right)^3 \right)^T, \tag{5.3}$$

and as in Section 2, $G(t) = V^{-1} G_* V$, $F(t) = V^{-1} F_*(t)$.

Finally, we set in equation (2.8), $G_- = G - \text{Im}(diag(G))$ and $L(t, z) = \|V^{-1}\| (\mu_1 \delta_1 + \mu_2 \delta_2) z^3$, where,

$$\delta_1 = \mu_1 \left( \sum_{k=1}^{4} abs(v_{2k}) \right)^3, \quad \delta_2 = \mu_2 \left( \sum_{k=1}^{4} abs(v_{4k}) \right)^3, \tag{5.4}$$

and, as in Section 1.2, $abs(q) = \sqrt{a^2 + b^2}$ if $q = a + ib$, and $abs(q) = |q|$ if $\text{Im}(q) = 0$.

Let us show briefly how to derive listed above formular for $L(t, z)$. Assume that $Y_i = \sum_{k=1}^{4} v_{ik} y_k$, $i = 2, 4$. Then (5.3) can be written as follow, $f = V^{-1} Y$, where $Y = (0 \quad -\mu_1 Y_2^3 \quad 0 \quad -\mu_2 Y_4^3)^T$. Subsequently, application of standard inequalities, $\|V^{-1} Y\| \leq \|V^{-1}\| \|Y\|$, $\|Y\|_2 \leq \|Y\|_1 = abs(\mu_1 Y_2^3) + abs(\mu_1 Y_4^3)$ and $\|y_k^m\| \leq \|y\|^m$, $m \in \mathbb{N}$, returns that $\|Y\|_1 = abs\left(\mu_1 \left(\sum_{k=1}^{4} v_{2k} y_k\right)^3\right) + abs\left(\mu_2 \left(\sum_{k=1}^{4} v_{4k} y_k\right)^3\right) = (\mu_1 \delta_2 + \mu_2 \delta_4) \|y\|^3$.

Consequently, $\|V^{-1} Y\| \leq \|V^{-1}\| (\mu_1 \delta_2 + \mu_2 \delta_4) z^3$.

Fourth order Runge-Kotter method with variable step size was used in our simulations. To estimate the boundaries of the trapping/ stability regions, the initial fourth-order systems were written in double polar – like coordinates as follows

$$x_1(t, x_0) = r_1(t, r_0, \varphi_0) \cos(\varphi_1(t, r_0, \varphi_0)), \quad x_2(t, x_0) = r_1(t, r_0, \varphi_0) \sin(\varphi_1(t, r_0, \varphi_0)),$$
$$x_3(t, x_0) = r_2(t, r_0, \varphi_0) \cos(\varphi_2(t, r_0, \varphi_0)), \quad x_4(t, x_0) = r_2(t, r_0, \varphi_0) \cos(\varphi_2(t, r_0, \varphi_0)),$$
$$r_0 = [r_{1,0} \quad r_{2,0}]^T, \quad \varphi_0 = [\varphi_{1,0} \quad \varphi_{2,0}]^T$$

To reduce the scope of simulations, we present below only computations of projections of the boundedness/stability regions on some coordinate planes. In our simulations, both angle coordinates were discretized with step equal $\pi/60$. For each fixed set of angle coordinates, the radial ones were adjusted sequentially to approximate the values located on the region's boundary. We started from small radial values and advanced them until rapid increases in $\|x(t, x_0)\|$ in two consecutive time - steps were numerically detected. Similar approach was used to approximate the threshold values in simulations of our scalar auxiliary equations. The running time in the latter computations practically does not depend upon dimension of the original system and is reduced further since $\|z(t, z_0)\|$ monotonically increases in $z_0$ for $\forall t \geq t_0$. In contrast, the running time required to estimate the boundary of trapping/stability regions increases with dimensions as $m^n$, where $m$ is the number of points that were taken to discretize a phase-space variable and $n$ is the number of these coordinates.

The results of simulation of the Van der – Pol -like system with the indicated parameters are partly shown in Fig. 1. Fig. 1a contrasts in solid, dashed, and dotted lines time-histories of the norms of solutions to equations (1.2), (2.10) and (2.7) with $F_0 = 0$, respectively, which correspond to a coupled Van der- Pol-like system. Note that in this case $x_0$ is located in the main part of the actual stability region and afar from its boundary. Clearly, the utility of (2.10) notably enhances the estimates delivered through the application of a more conservative equation (2.7) with $F_0 = 0$ and adequatelly estimates the norm of actual solutions to (1.2). In the subsequent simulations presented in this subsection we set that, $\mu_1 = \mu_2 = 18.2$. Figs. 1b - d display in dashed and solid lines the projections on some coordinate planes of the boundaries of stability regions of trivial solutions to equations (1.2) and (2.10) corresponding to the Van der Pol- like system. In fact, simulations of (2.10) estimate the thereshold value of $z_0$, i.e, $\bar{\zeta}$ which, in turn, defines the ellipsoid imbedded in the stability region through application of the formula, $\bar{z} = \|V^{-1} x_0\|$. Clearly, the attained estimates of the boundary of stability region turn out to be fairly concervative if the structure and parameters of the model are defined presicelly. Yet, the practical merit of these estimates becomes more apparent for systems under uncertainty, where, for instance, it is assumed that at most the norms of some model components or parameters are defined presizelly.

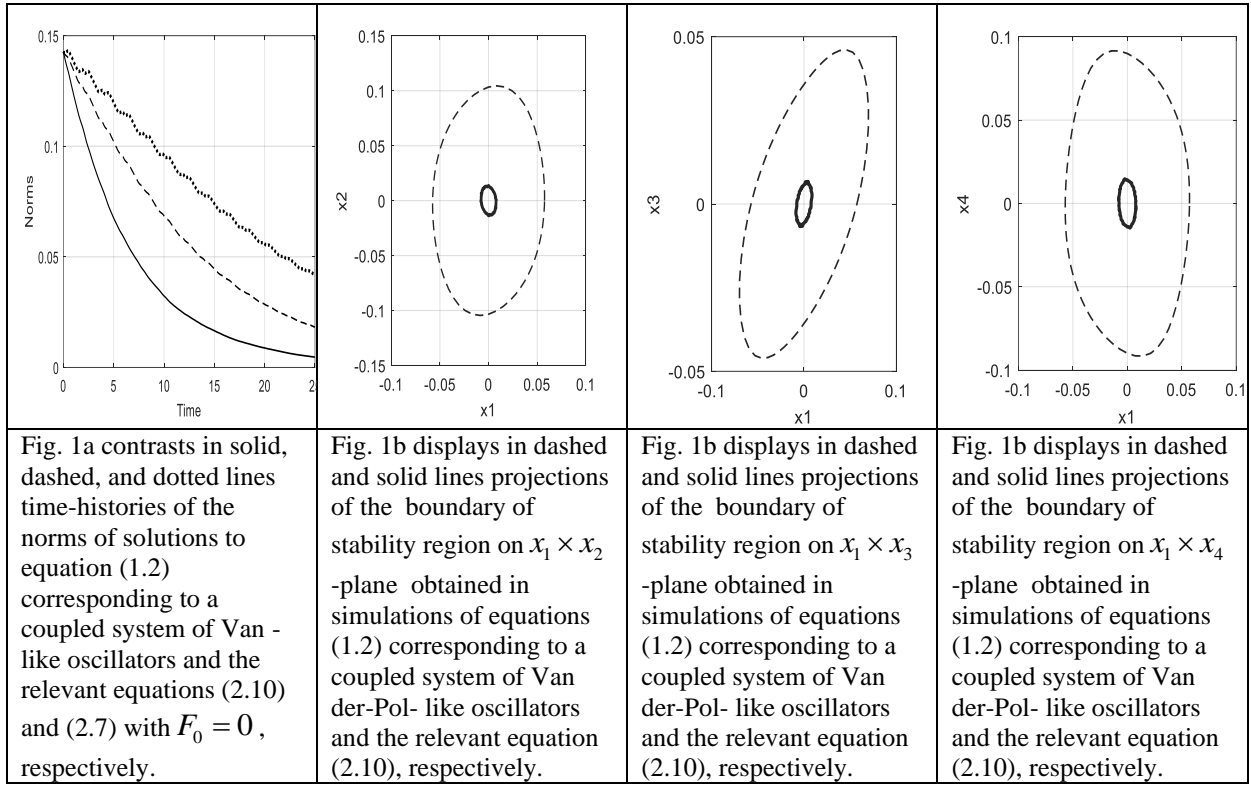

| Fig. 1a contrasts in solid, dashed, and dotted lines time-histories of the norms of solutions to equation (1.2) corresponding to a coupled system of Van - like oscillators and the relevant equations (2.10) and (2.7) with $F_0 = 0$, respectively. | Fig. 1b displays in dashed and solid lines projections of the boundary of stability region on $x_1 \times x_2$ -plane obtained in simulations of equations (1.2) corresponding to a coupled system of Van der-Pol- like oscillators and the relevant equation (2.10), respectively. | Fig. 1b displays in dashed and solid lines projections of the boundary of stability region on $x_1 \times x_3$ -plane obtained in simulations of equations (1.2) corresponding to a coupled system of Van der-Pol- like oscillators and the relevant equation (2.10), respectively. | Fig. 1b displays in dashed and solid lines projections of the boundary of stability region on $x_1 \times x_4$ -plane obtained in simulations of equations (1.2) corresponding to a coupled system of Van der-Pol- like oscillators and the relevant equation (2.10), respectively. |
|---|---|---|---|

## 5.2 Coupled System of Duffing -like Oscillators

For this system, (5.1) remains intact, but (5.2) should be altered as follows,

$$f_\times = \begin{pmatrix} 0 & -\mu_1 x_1^3 & 0 & -\mu_2 x_3^3 \end{pmatrix}^T \tag{5.5}$$

which, in turn, requires some alteration in (5.4), i.e., $v_{2k} \to v_{1k}$ and $v_{4k} \to v_{3k}$, respectively. Fig. 2 displays in dashed and solid lines the projections on some coordinate -planes of the boundaries of stability regions of the trivial solutions to equations (1.2) and (2.8) corresponding to a coupled system of Duffing – like oscillators, respectively. As in Section 5.1, the current estimates turn out to be rather conservative for precisely defined systems but appear to be more compelling for systems under uncertainties. Note that the simulations on Fig. 2b are started at $t_0 = 1.43$ whereis other plots in this figure are simulated with $t_0 = 0$. Clearly, comparison of Figures 2a and 2b shows that the shapes of the attractors' boundaries in these figures are somewhat affected by the initial time moment. However, such dependence might become less pronouced for other values of $t_0$ since time-dependent coefficients are almost periodic in this model.

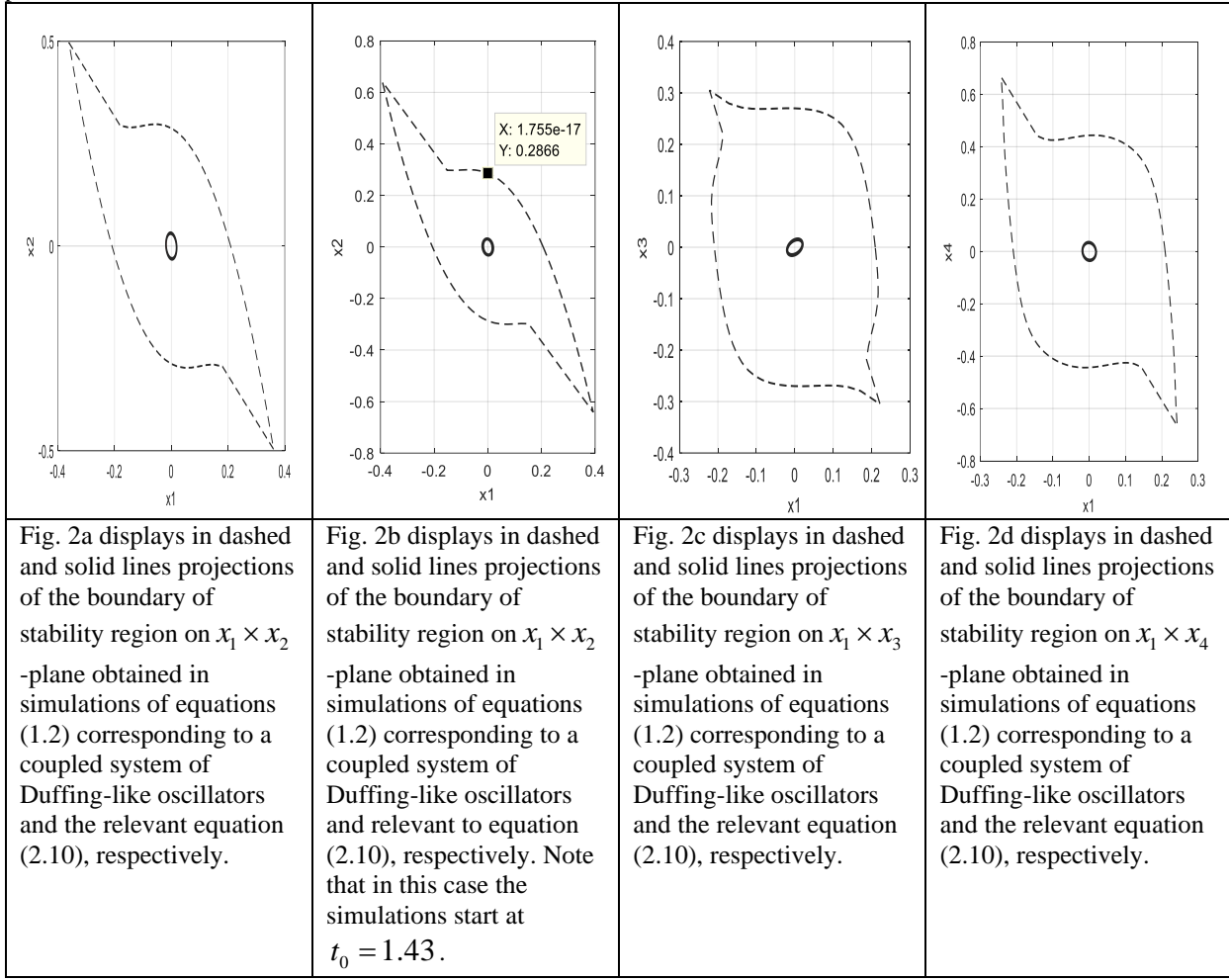

| Fig. 2a displays in dashed and solid lines projections of the boundary of stability region on $x_1 \times x_2$ -plane obtained in simulations of equations (1.2) corresponding to a coupled system of Duffing-like oscillators and the relevant equation (2.10), respectively. | Fig. 2b displays in dashed and solid lines projections of the boundary of stability region on $x_1 \times x_2$ -plane obtained in simulations of equations (1.2) corresponding to a coupled system of Duffing-like oscillators and relevant to equation (2.10), respectively. Note that in this case the simulations start at $t_0 = 1.43$. | Fig. 2c displays in dashed and solid lines projections of the boundary of stability region on $x_1 \times x_3$ -plane obtained in simulations of equations (1.2) corresponding to a coupled system of Duffing-like oscillators and the relevant equation (2.10), respectively. | Fig. 2d displays in dashed and solid lines projections of the boundary of stability region on $x_1 \times x_4$ -plane obtained in simulations of equations (1.2) corresponding to a coupled system of Duffing-like oscillators and the relevant equation (2.10), respectively. |
|---|---|---|---|

Next, we review some results of simulations on estimating the trapping regions for nonhomogeneous equation (1.1). Fig. 3 displays in dashed and solid lines projections on the coordinate planes of the boundaries of trapping regions about zero for equations (1.1) and (2.8) that are computed for a coupled Duffing – like system with variable coefficients. In these simulations, $\mu_1 = \mu_2 = 0.3$ and $F_{01} = 0.01$, $F_{02} = 0$. Our outcomes in these simulations are comparable to the ones that were presented prior and show that our fairly conservative estimates of the boundary of trapping region of (1.1) may become more practically appealing for a system under uncertainty. Yet, the utility of these estimates for precisely defined systems provides a numerical affirmation of our boundedness and stability criteria. Nonetheless, combining the techniques outlined in this paper with the pertinent successive approximations should considerably refine the accuracy of the corresponding estimates, see [14], where such a strategy was developed and implemented for suitable systems.

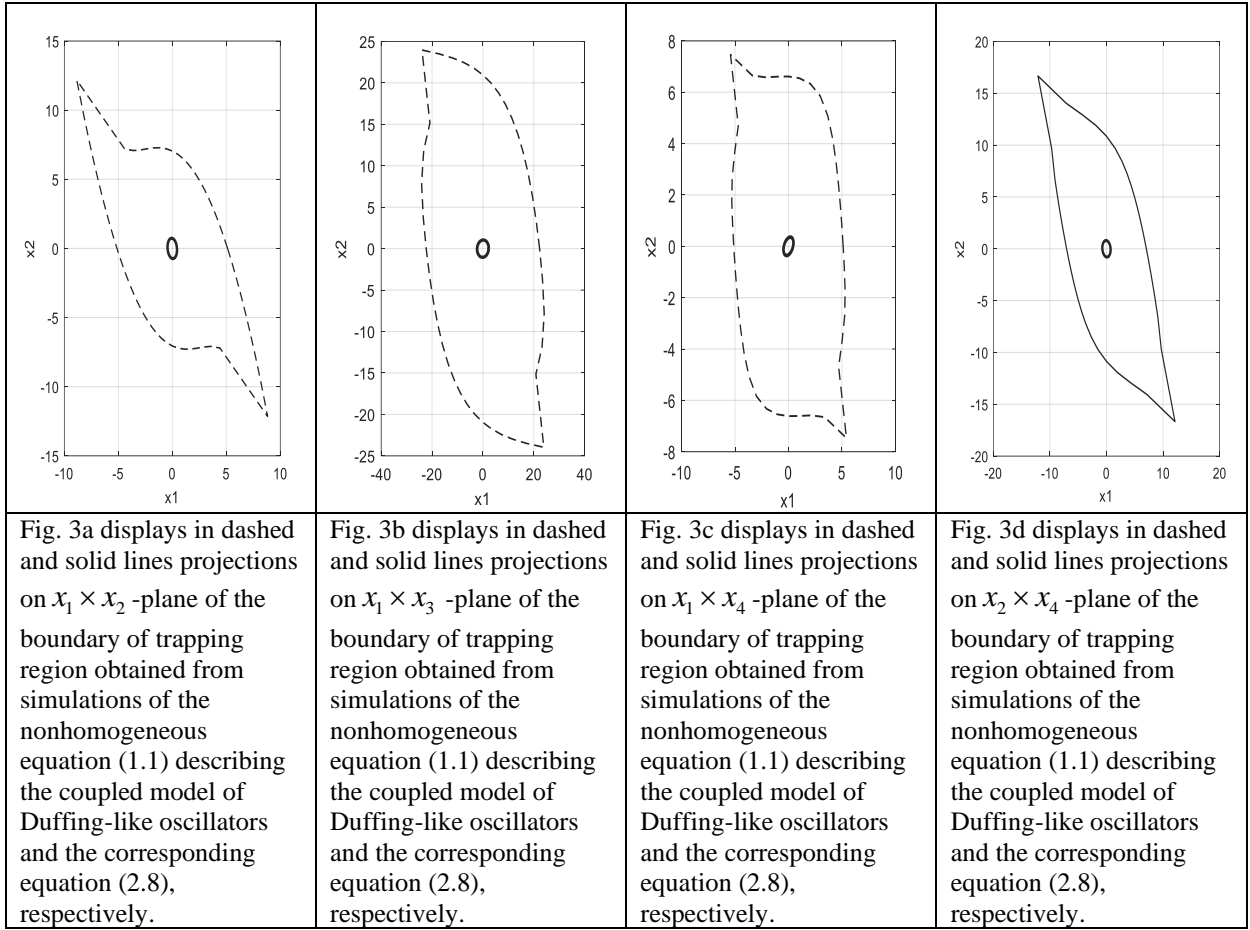

| Fig. 3a displays in dashed and solid lines projections on $x_1 \times x_2$ -plane of the boundary of trapping region obtained from simulations of the nonhomogeneous equation (1.1) describing the coupled model of Duffing-like oscillators and the corresponding equation (2.8), respectively. | Fig. 3b displays in dashed and solid lines projections on $x_1 \times x_3$ -plane of the boundary of trapping region obtained from simulations of the nonhomogeneous equation (1.1) describing the coupled model of Duffing-like oscillators and the corresponding equation (2.8), respectively. | Fig. 3c displays in dashed and solid lines projections on $x_1 \times x_4$ -plane of the boundary of trapping region obtained from simulations of the nonhomogeneous equation (1.1) describing the coupled model of Duffing-like oscillators and the corresponding equation (2.8), respectively. | Fig. 3d displays in dashed and solid lines projections on $x_2 \times x_4$ -plane of the boundary of trapping region obtained from simulations of the nonhomogeneous equation (1.1) describing the coupled model of Duffing-like oscillators and the corresponding equation (2.8), respectively. |
|---|---|---|---|

Comparison of our latter [13] and current techniques shows that the current one gains computational efficiency and applicable to an essentially wider class of systems. Nonetheless, the accuracy of this technique decreases if $|\alpha_1|$ is a small value, $\sup_{\forall t \geq t_0} \|G_-(t)\|$ is a large value, and $\sup_{\forall t \geq t_0} L(t,z)$ becomes sufficiently large for relevant values of $z$. The accuracy of our prior technique is practically non sensitive to variations of $\sup_{\forall t \geq t_0} \|G_-(t)\|$ but decreases with growth of $\sup_{\forall t \geq t_0} L(t,z)$ as well. However, this prior technique also requires numerical simulation of $w(t)$, $p(t)$ and $c(t)$ which is avoided in the current approach.

## 6. Conclusion and Forthcoming Studies

In [13], we derived a scalar nonlinear auxiliary equation with solutions bounded from above the norm of solutions to some nonlinear and nonautonomous systems and applied this technique to the assessment of the boundedness and stability of solutions to such systems. Yet, the application domain of this approach is constrained since the underlying auxiliary equation contains a function that in some cases can approach infinity with $t \to \infty$. The current paper presents a novel technique casting the auxiliary equation in a modified form escaping this limitation for a wide class of systems arising in applications. Furthermore, the modified auxiliary equation is simpler than its former counterpart and more efficient in computation. Still, in the common application domain, the current estimates frequently turn out to be more conservative than their prior analogs.

Next, we present various novel boundedness and stability criteria which are stemmed from analysis of the modified nonlinear auxiliary equation and its linear counterpart. The latter approach leads to close-form and simplified boundedness and stability criteria that should be valuable in resolving relevant control problems.

Lastly, we authenticate our study using inclusive simulations of systems with typical dissipative and conservative nonlinear components that were intractable to our prior technique. The simulations show that our upper estimates of the norms of solutions to nonlinear and nonautonomous systems turn out to be adequate if they stem from the central parts of the corresponding trapping/stability regions. Nonetheless, the boundaries of trapping/ stability regions are estimated rather conservatively by our current technique if the corresponding systems are defined precisely. Yet, such estimates turn out to be more appealing in applications for systems under uncertainty.

Note also that the precision of our estimates can be substantially improved through the integration of our current technique with a methodology of successive approximations that was outlined in [14]. This should provide bilateral bounds for the norms of solutions to a broad class of time – varying and nonlinear systems and enhance recursively the accuracy of estimation of their trapping/ stability regions.

**Appendix**

Assume that $M(\lambda) = \lambda + \max_k |\alpha_k - \lambda|$ and show that $\min_\lambda (M(\lambda))$ is attained at $\lambda = \alpha_n$ under our prior assumption that $\alpha_i > \alpha_{i+1}$, $i = 1,...,n-1$. For this sake, we assume in order that:

1. $\lambda \geq \alpha_1$, then, $|\alpha_i - \lambda| = -(\alpha_i - \lambda)$, $i = 1,...,n$ which yields that $\min_\lambda (M(\lambda)) = (2\lambda - \alpha_i)_{\lambda = \alpha_1} \geq \alpha_1$

2. $\lambda \leq \alpha_n$, then, $|\alpha_i - \lambda| = (\alpha_i - \lambda)$, $i = 1,...,n$ which yields that $\min_\lambda (M(\lambda)) = \alpha_n$.

3. $\alpha_{k+1} \leq \lambda \leq \alpha_k$, $k = 1,...,n-1$. This assertion splits into the following cases:

3a. $i \geq k+1$, then $|\alpha_i - \lambda| = (\alpha_i - \lambda)$ which yields that $\min_\lambda (M(\lambda)) = \alpha_{k+1}$.

3b. $i < k+1$, then $|\alpha_i - \lambda| = -(\alpha_i - \lambda)$, $i = 1,...,n$ which yields that
$\min_\lambda (M(\lambda)) = (2\lambda - \alpha_i)_{\lambda = \alpha_k} = 2\alpha_k - \alpha_i \geq \alpha_k$.

Lastly, we inference that $\min_\lambda (M(\lambda)) = \alpha_n$ and $M(\alpha_n) = \alpha_n + \max_k |\alpha_k - \alpha_n| = \alpha_1$ □

**Acknowledgement**. The programs used in the simulations provided in Section 5 of this paper were developed by Steve Koblik.

Data sharing is not applicable to this article as no datasets were generated or analyzed during the current study.

This research was not supported by external funds or grants.

The author declares that he has no conflict of interest.